\documentclass[article]{amsart}
\usepackage{cases}
\usepackage{amsmath, amsfonts, pifont, amssymb}
\usepackage{verbatim}
\usepackage{geometry}
\newtheorem{theorem}{Theorem}[section]
\newtheorem{lemma}[theorem]{Lemma}
\newtheorem{proposition}[theorem]{Proposition}
\newtheorem{corollary}[theorem]{Corollary}

\newcommand{\R}{\mathbb{R}}

\newcommand{\f}{\frac}
\newcommand{\beq}{\begin{equation}}
\newcommand{\eeq}{\end{equation}}
\newcommand{\beqq}{\begin{equation*}}
\newcommand{\eeqq}{\end{equation*}}

\theoremstyle{definition}
\newtheorem{definition}[theorem]{Definition}

\theoremstyle{remark}

\newtheorem{conjecture}[theorem]{Conjecture}
\numberwithin{equation}{section}

\numberwithin{equation}{section}

\begin{document}

\title{Dynamics of subcritical threshold solutions for energy-critical NLS }

\author{Qingtang Su and Zehua Zhao}

\maketitle

\begin{abstract}
In this paper, we study the dynamics of subcritical threshold solutions for focusing energy critical NLS on $\mathbb{R}^d$ ($d\geq 5$) with nonradial data. This problem with radial assumption was studied by T. Duyckaerts and F. Merle in \cite{DM} for $d=3,4,5$ and later by D. Li and X. Zhang in \cite{LZ} for $d \geq 6$. We generalize the conclusion for the subcritical threshold solutions by removing the radial assumption for $d\geq 5$. A key step is to show exponential convergence to the ground state $W(x)$ up to symmetries if the scattering phenomenon does not occur. Remarkably, an interaction Morawetz-type estimate are applied. 
\end{abstract}
\bigskip

\bigskip

\section{Introduction}
We consider the following focusing energy-critical nonlinear Schr{\"o}dinger initial value problem:
\begin{equation}\label{maineq}
\aligned
(i\partial_t+ \Delta_{\mathbb{R}^d}) u &= F(u) = \lambda|u|^{p} u, \quad (t,x)\in I\times \mathbb{R}^d\\
u(0,x) &= u_{0} \in \dot{H}^{1}(\mathbb{R}^d), \textmd{ where } p=\frac{4}{d-2},d\geq 5,\lambda=-1.
\endaligned
\end{equation}
If $\lambda=1$, the equation would be defocusing (see \cite{Iteam1,ec4d,V1} for some results regarding defocusing case). Cauchy problem \eqref{maineq} has been studied in \cite{CW}. Locally, there exists a unique solution defined on a maximal interval $I$ such that for strictly smaller subinterval $J$ of $I$,

\begin{equation}
||u||_{S(J)}<\infty, 
\end{equation}
where $S(J):=L^{\frac{2(d+2)}{d-2}}(J\times \mathbb{R}^d)$ (scattering norm). Moreover, the initial value problem (1.1) is called energy-critical problem since the energy of the solution is invariant under the scaling symmetry in the following sense. First, the energy of the solution
\begin{equation}
E(u(t))=\frac{1}{2} \int |\nabla u(t,x)|^2 dx-\frac{1}{2^{*}}\int |u(t,x)|^{2^{*}} dx \textmd{ where } 2^{*}=\frac{2d}{d-2}
\end{equation}
is a conserved quantity. Furthermore, the solutions of (1.1) are invariant under the following transformation: for a solution $u(t,x)$,
\begin{equation}
\frac{e^{i\theta_0}}{\lambda_0^{(d-2)/2}} u(\frac{t_0+t}{\lambda_0^2},\frac{x_0+x}{\lambda_0}) \textmd{ where } (\theta_0,\lambda_0,t_0,x_0)\in (\mathbb{R}\times (0,\infty) \times \mathbb{R}\times \mathbb{R}^d) 
\end{equation}
is also a solution. The transformation group is generated by translations, rotations and scalings according to the symmetric structure of the equation (1.1). It is straightforward to verify that these transformations preserve the $S(\mathbb{R})$-norm, as well as the $\dot{H}^1$-norm, the $L^{2^*}$-norm and thus the energy.\vspace{2mm} 

Generally speaking, (1.1) is a special case (when critical index $s_c=1$) of the following critical initial value problem:
\begin{equation}\label{generaleq}
\aligned
(i\partial_t+ \Delta) u &= F(u) = -|u|^{p} u, \quad (t,x)\in I\times \mathbb{R}^d\\
u(0,x) &= u_{0} \in \dot{H}^{s_c}(\mathbb{R}^d), \textmd{ where } p:=\frac{4}{d-2s_c}.
\endaligned
\end{equation}
For focusing energy-critical NLS, there is an important radial stationary solution $W$ (see \cite{Au,Ta} for more information) satisfying the following elliptic equation:
\begin{equation}
\Delta_{\mathbb{R}^d} W=-|W|^{\frac{4}{d-2}}W,
\end{equation}
\noindent and the explicit expression of $W$ is
\begin{equation}
W(x)=\frac{1}{(1+\frac{|x|^2}{d(d-2)})^{\frac{d-2}{2}}}.
\end{equation}
We are interested in the dynamics of the solutions to Cauchy problem (1.1). First, for the defocusing case (when $\lambda=1$ in \eqref{maineq}), there is a scattering result as follows:

\begin{theorem}[Scattering for defocusing energy-critical NLS]\label{defocmain}
For $d\geq 3$, assume $u_0 \in \dot{H^{1}}(\mathbb{R}^d)$, then there exists a unique global solution $u \in C(\mathbb{R}:\dot{H^{1}}(\mathbb{R}^4))$ of the initial-value problem 
\begin{equation}
(i\partial_{t} +\Delta_{\mathbb{R}^d})u=u|u|^{\frac{4}{d-2}}, \quad u(0)=u_0.
\end{equation}

\noindent Moreover, this solution scatters in the sense that there exists $\psi^{\pm \infty} \in \dot{H^{1}}(\mathbb{R}^d)$ such that 
\begin{equation}\label{euqation}
\|u(t)-e^{it\Delta} \psi^{\pm \infty}\|_{\dot{H^{1}}(\mathbb{R}^d)} \to 0
\end{equation}
as $t \to \pm \infty$. 
\end{theorem}
\emph{Remark. }Theorem \ref{defocmain} is proved in \cite{Iteam2} for $d=3$, in \cite{ec4d} for $d=4$ and in \cite{V1} for $d\geq 5$.\vspace{2mm}

However, different from the defocusing case, the dynamics of the solutions of focusing energy-critical NLS are much richer and tightly dependent on the sizes of initial data. Compared with the ground state $W$, we may roughly consider three scenarios, i.e. the initial data is `below', `at' and `above' the ground state $W$ in the sense of energy.\vspace{2mm} 

First, if the initial data is `below' the ground state in the sense of $E(u_0)<E(W)$, there is a famous conjecture as follows:
\begin{conjecture}[Ground state conjecture for energy-critical NLS]\label{gsc}
When $ d\geq 3$, we consider the initial value problem \eqref{maineq}, assuming initial data $u_0 \in \dot{H}^1(\mathbb{R}^d)$, under the assumption that the solution $u(t)$ satisfies
\begin{equation}\label{asp1}
\sup\limits_{t \in I}||u(t)||_{\dot{H}^1(\mathbb{R}^d)}<||W||_{\dot{H}^1(\mathbb{R}^d)},
\end{equation}
where $I$ is the lifespan of the solution $u(t)$, then there exists a unique global solution $u \in C(\mathbb{R}:\dot{H}^1(\mathbb{R}^d))$ of the energy critical initial-value problem,
\begin{equation}
(i\partial_t+\Delta_{\mathbb{R}^d})u=-u|u|^{\frac{4}{d-2}}, \quad u(0)=u_0.
\end{equation} 
Moreover, this solution scatters in the sense that there exists $\phi^{\pm} \in \dot{H}^1(\mathbb{R}^d)$, such that
\begin{equation}
||u(t)-e^{it\Delta_{\mathbb{R}^d}}\phi^{\pm}||_{\dot{H}^1} \rightarrow 0,\quad \textmd{as} \quad t\rightarrow \pm \infty.
\end{equation}

\end{conjecture}
\emph{Remark.} In Conjecture \ref{gsc}, according to energy trapping theorem (see \cite{KM1}), we can replace the priori assumption (\ref{asp1}) by the following assumption regarding the initial data,
\begin{equation}
||u_0||_{\dot{H}^1}<||u_0||_{\dot{W}^1}, \quad E(u_0)<E(W). 
\end{equation}

\emph{Remark.} The main idea of Conjecture \ref{gsc} is, if the initial data is below the ground state, then the dynamics of the solutions would be similar to the defocusing case, i.e. the behavior of the solutions resemble linear solutions. Conjecture \ref{gsc} was proved with radial assumption by C. Kenig and F. Merle for $d=3,4,5$ (see \cite{KM1}). Later, Conjecture \ref{gsc} was proved for $d\geq 5$ by R. Killip and M. Visan (see \cite{KV ecg}) and for $d=4$ by B. Dodson (see \cite{BD ecg}). The case when $d=3$ is still open.\vspace{2mm} 

\emph{Remark.} We refer to \cite{BD mcg} for the analogue of Conjecture \ref{gsc} regarding the mass-critical case.\vspace{2mm}

Also, we are interested in the situation when the solutions are `at' the ground state in the sense of $E(u_0)=E(W)$. We call these solutions ``energy threshold solutions''. Our ultimate goal is to give a classification of the threshold solutions of \eqref{maineq} with critical energy $E(W)$. The following theorem is the main theorem of \cite{DM} and \cite{LZ}.
\begin{theorem}\label{DMLZ}
\noindent For $d \geq 3$, let $u$ be the solution of (1.1) with radial initial data $u_0$ satisfying 
\begin{equation}
E(u_0)=E(W)=\frac{1}{dC_d^d}
\end{equation}
and $I$ its maximal interval of definition. Then the following conclusions hold:\vspace{2mm}

 (a) If $\int |\nabla u_0|^2<\int |\nabla W|^2=\frac{1}{dC_d^d}$ then $I=\mathbb{R}$. Moreover, either $u=W^{-}$ up to the symmetry of the equation, or $u$ scatters in both time directions.\vspace{2mm}

 (b) If $\int |\nabla u_0|^2=\int |\nabla W|^2$ then $u=W$ up to symmetry of the equation.\vspace{2mm}

 (c) If $\int |\nabla u_0|^2>\int |\nabla W|^2$ and $u_0\in L^2$ then either $u=W^{+}$ up to symmetry of the equation or $I$ is finite.
\end{theorem}

\emph{Remark.} Theorem \ref{DMLZ} is proved for $d=3,4,5$ in \cite{DM} (by T. Duyckaerts and F. Merle) and in \cite{LZ}  for $d\geq 6$ (by D. Li and X. Zhang). \vspace{2mm}

\emph{Remark.} $C_d$ is the best constant in the Sobolev inequality for $d$-dimensional case (see \cite{Au,Ta}). \vspace{2mm}

\emph{Remark.} We refer to \cite{DM2,LZ2} for the analogues of this result for nonlinear wave equations.\vspace{2mm}

 \emph{Remark.} In Theorem \ref{DMLZ}, stationary, radial functions $W^{-}$ and $W^{+}$ satisfy following properties (see \cite{DM,LZ} for more information):
 \begin{theorem}
 For $d\geq 3$, there exists radial solutions $W^{-}$ and $W^{+}$ of \eqref{maineq} such that
 \begin{equation}
 E(W)=E(W^{-})=E(W^{+}),
 \end{equation}
  \begin{equation}
 T_{+}(W^{-})=T_{+}(W^{+})=+\infty,\quad W^{\pm}(t)=W \textmd{ in }\dot{H}^1,
 \end{equation}
  \begin{equation}
  ||W^{-}||_{\dot{H}^1}<||W||_{\dot{H}^1},T_{-}(W^{-})=+\infty,||W^{-}||_{S((-\infty,0])}<+\infty,
 \end{equation}
   \begin{equation}
 ||W^{+}||_{\dot{H}^1}>||W||_{\dot{H}^1}.
 \end{equation}
 \end{theorem}
\emph{Remark.} In this paper, we use the same $W^-$. Regarding the construction of $W^-$, we refer to \cite{DM} for $d=5$ and \cite{LZ} for $d\geq 6$.  \vspace{2mm}

It is natural to consider the nonradial case by removing the radial assumption in Theorem \ref{DMLZ}. When the energy of the solution equals the energy of the ground state $W$, as discussed above, there are three cases dependent on the kinetic energy ($\dot{H}^1$-norm) of the initial data.\vspace{2mm} 

First, for case (b) (kinetic energy critical), similar conclusion still holds for the nonradial case. We recall the following result (see \cite{Au,Ta}): 
\begin{theorem}\label{sobolev}
Let $C(d)$ denote the sharp constant in Sobolev inequality,
\begin{equation}
||f(x)||_{L^{\frac{2d}{d-2}}_x(\mathbb{R}^d)}\leq C(d)||\nabla f||_{L^2(\mathbb{R}^d)}.
\end{equation}
Then the equality holds if and only if $f=W$ up to symmetries in the following sense,
\begin{equation}
f(x)=e^{i\theta_0}\lambda_0^{-\frac{d-2}{2}}W(\frac{x-x_0}{\lambda_0})
\end{equation}
for $(\theta_0,\lambda_0,x_0) \in \mathbb{R}\times \mathbb{R}^{+} \times \mathbb{R}^d$.
Thus, in particular, if $u$ is a solution of \eqref{maineq} satisfying 
\begin{equation}
E(u_0)=E(W),\quad \int |\nabla u_0|^2=\int |\nabla W|^2=\frac{1}{dC_d^d}.
\end{equation}
Then $u_0$ coincides with $W$ up to symmetries, so does the corresponding solution $u(t)$. 
\end{theorem}

For the other two cases, the conclusion is nontrivial. The full resolution of case (c) (kinetic supercritical case) seems to require some new techniques and it is very different from case (a) (kinetic subcritical case), so we leave it for a future work. In this paper, we consider the dynamics of subcritical threshold solutions (case (a)) and the main result is as follows:
\begin{theorem}\label{main}
When $d\geq 5$, let $u$ be the solution of (1.1) with initial data $u_0$ satisfying 
\begin{equation}\label{subcr}
E(u_0)=E(W)=\frac{1}{dC_d^d}, \quad  \int |\nabla u_0|^2<\int |\nabla W|^2=\frac{1}{dC_d^d},
\end{equation}
and $I$ its maximal interval of definition. Then $u$ is global, i.e. $I=\mathbb{R}$. Moreover, either $u=W^{-}$ up to the symmetry of the equation, or $u$ scatters in both time directions.
\end{theorem}
\emph{Remark. }Theorem \ref{main} gives a classification of the subcritical threshold solutions. As for the dynamics of the subcritical threshold solutions, there are exact two situations: the solution scatters in two directions or the solution equals to $W^-$ up to symmetries. \vspace{2mm}

The road map of approaching Theorem \ref{main} is briefly explained as follows. First, we show the subcritical threshold solution is global. Moreover, if the solution does not scatter, the solution is almost periodic in the sense of satisfying the compactness condition. At last, we show the exponential convergence to the ground state $W(x)$ and use it to obtain the main theorem.\vspace{2mm}

One main difference from the radial case is the appearance of the translation parameter $x(t)$ in the nonradial setting. We need to deal with $x(t)$ carefully which arises from the compactness argument (Theorem \ref{compactness}). For the radial case, the translation parameter is trivially $0$. This difference causes changes in subsequent arguments since we need to control the translation parameter. Thus we need to establish a compactness result (Theorem \ref{compactness}) and an orthogonal decomposition (Theorem 4.5) in the nonradial setting and apply an interaction Morawetz estimate to obtain the exponential convergence (Theorem \ref{4main}).\vspace{2mm}

At last, we refer to \cite{KM1} for the situation when the solutions are `above' the ground state in the sense of $E(u_0)>E(W)$.\vspace{2mm}

\textbf{Organization of this paper:} In Section 1, we introduce the background, existing results and the main result of this paper; in Section 2, we discuss preliminaries, basic tools and the compactness result; in Section 3, we prove that the scaling function $\lambda(t)$ in the compactness argument has a lower bound and use it to obtain some important properties regarding the almost periodic solution; in Section 4, we prove the exponential convergence to ground state $W$ for subcritical threshold solutions if the scattering phenomenon does not occur; in Section 5, we use the exponential convergence result established in Section 4 and the results in \cite{DM,LZ} to prove the main theorem; in Section 6 (Appendix), we give the proofs of Lemma 4.5 and Lemma 4.6.
 
\section{Preliminaries and compactness result }
In this section, we discuss preliminaries, basic tools and the compactness result (Theorem \ref{compactness}).\vspace{2mm}

We write $X\lesssim Y$  or $Y\gtrsim  X$ whenever $X\leq CY$  for some constant $C>0$. Moreover, we use $O(Y)$ to denote any quantity $X$ such that  $|X|\lesssim Y$ and use $x(t)=o(t)$ to denote a time-dependent quantity $x(t)$ such that $x(t)/t \rightarrow 0$ as $t\rightarrow \infty$. If $X\lesssim Y$  and $Y \lesssim X$  hold simultaneously, we abbreviate that by $X\sim Y.$  Without special clarification, the implicit constant $C$  can vary from line to line. We use Japanese bracket $\langle x\rangle$ to denote $(1+|x|^2)^{\f{1}{2}}.$ \vspace{2mm}
 
Define the Fourier transform on $\R^d$ by
\begin{align*}
  \hat{f(\xi)}:=(2\pi)^{-\f{d}{2}}\int_{\R^d}e^{-ix\xi}f(x)dx,
\end{align*}
and the homogeneous Sobolev norm as
\begin{align*}
  \|f\|_{\dot{H}^s(\R^d)}:=\||\nabla|^s f\|_{L_x^2(\R^d)}
\end{align*}
where
\begin{align*}
  \widehat{|\nabla|^s f}(\xi):=|\xi|^s \hat{f}(\xi).
\end{align*}
Now we recall Littlewood-Pelay theory which is an important tool in the area of partial differential equations. Let $\phi(\xi)$  be a radial bump function supported in the ball $\{\xi\in \R^d :|\xi|\leq \f{11}{10}\}$  and equals $1$  on the ball $\{\xi\in \R^d: |\xi|\leq 1\}.$  For each dyadic number $N>0,$  we define
\begin{align*}
    \widehat{P_{\leq N}f}(\xi):=&\varphi\big(\frac{\xi}{N}\big)\hat{f}(\xi), \\
    \widehat{P_{> N}f}(\xi):=&\big(1-\varphi(\frac{\xi}{N})\big)\hat{f}(\xi), \\
    \widehat{P_{N}f}(\xi):=&\big(\varphi(\frac{\xi}{N})-\varphi(\frac{2\xi}{N})\big)\hat{f}(\xi),
\end{align*}
with similar definitions for $P_{< N}$ and $P_{\geq N}$. Also,  we define
\begin{align*}
    P_{M<\cdot\leq N}:=P_{\leq N}-P_{\leq M},
\end{align*}
whenever $M<N$. We state two useful results regarding the Littlewood-Paley operators as follows:
\begin{lemma}[Bernstein's inequalities]\label{2Bineq} For $1\leq r\leq q \leq \infty$, $s \geq 0$, we have
\begin{equation}
|||\nabla|^{\pm s} P_N f||_{L^r(\mathbb{R}^d)} \sim N^{\pm s}||P_N f||_{L^r(\mathbb{R}^d)},
\end{equation}
\begin{equation}
|||\nabla|^{ s} P_{\leq N} f||_{L^r(\mathbb{R}^d)} \lesssim N^s ||P_{\leq N} f||_{L^r(\mathbb{R}^d)},
\end{equation}
\begin{equation}
|| P_{\geq N} f||_{L^r(\mathbb{R}^d)} \lesssim N^{-s} |||\nabla|^{s} P_{\geq N} f||_{L^r(\mathbb{R}^d)},
\end{equation}
\begin{equation}
|| P_{\leq N} f||_{L^q(\mathbb{R}^d)} \lesssim N^{\frac{d}{r}-\frac{d}{q}} || P_{\leq N} f||_{L^r(\mathbb{R}^d)}.
\end{equation}
\end{lemma}

\begin{lemma}[Littlewood-Pelay square function estimate]\label{2lp} For $1<r<\infty$,
\begin{equation}
||(\sum |P_N f(x)|^2)^{1/2}||_{L^r_x(\mathbb{R}^d)} \sim ||f||_{L^r_x(\mathbb{R}^d)},
\end{equation}

\begin{equation}
||(\sum N^{2s}|P_N f(x)|^2)^{1/2}||_{L^r_x(\mathbb{R}^d)} \sim || |\nabla|^s f||_{L^r_x(\mathbb{R}^d)} \textmd{ for any } s,
\end{equation}
\begin{equation}
||(\sum N^{2s}|P_{>N} f(x)|^2)^{1/2}||_{L^r_x(\mathbb{R}^d)} \sim ||  |\nabla|^s f||_{L^r_x(\mathbb{R}^d)} \textmd{ for any } s>0.
\end{equation}

\end{lemma}
Then we recall dispersive estimate, Strichartz estimate and fractional product rule as follows.
\begin{lemma}[Dispersive estimate]\label{22disp}
\begin{equation}\label{disp}
||e^{it\Delta} f ||_{L^{\infty}_x(\mathbb{R}^d)} \lesssim |t|^{-\frac{d}{2}}||f||_{L^1(\mathbb{R}^d)}.
\end{equation}
\end{lemma}
\noindent \emph{Remark.} When the dimension of the function is higher, the decay is faster. Moreover, if we interpolate \eqref{disp} with $||e^{it\Delta}f||_{L^2(\mathbb{R}^d)}=||f||_{L^2(\mathbb{R}^d)}$ (Plancherel formula), we can obtain
\begin{equation}
||e^{it\Delta} f ||_{L^{r}_x(\mathbb{R}^d)} \lesssim |t|^{-(\frac{d}{2}-\frac{d}{r})}||f||_{L^{r^{'}}(\mathbb{R}^d)}.
\end{equation}
\noindent where $2\leq r \leq \infty$, $\frac{1}{r}+\frac{1}{r^{'}}=1$ and $t \neq 0$.

\begin{definition}[Admissible pair]\label{admi} Let $d\geq 5,$ we call a pair of exponent $(q,r)$ admissible if
\begin{align}
  \frac{2}{q}=d(\f{1}{2}-\f{1}{r}) \quad \textrm{with}\quad 2\leq q\leq \infty.
\end{align}
For a time interval $I$, we define
\begin{align}
  \|u\|_{\rm {S}(I)}:=\sup\{\|u\|_{L^q_tL^r_x(I\times\mathbb{R}^d)}:(q,r)  \text{   admissible}\}.
\end{align}
We also define the dual of $ S(I)$ by $N(I)$. Note that
\begin{align}
  \|u\|_{\rm N(I)}\lesssim \|u\|_{L_t^{q'}L_x^{r'}(I\times \R^d)}\quad \text{for any admissible pair}~(q,r).
\end{align}
\end{definition}

\begin{lemma}[Strichartz estimate]
  Let ~$u:I\times \mathbb{R}^d\rightarrow \mathbb{C}$~be a solution to
  \begin{align}
    (i\partial_t+\Delta)u=F
  \end{align}
  and let~$ s\geq 0$, then
  \begin{align}
    \||\nabla|^su\|_{S(I)}\lesssim\|u(t_0)\|_{\dot{H}^s_x}+\||\nabla|^sF\|_{N(I)},
  \end{align}
  for any ~$t_0 \in I$.
\end{lemma}

\begin{lemma}[Fractional chain rule]\label{chain}
Suppose ~$G\in C^1(\mathbb{C})$~and ~$s\in (0,1].$~Let ~$1<r<r_2<\infty $~and ~$1<r_1\leq \infty $~be such that~$\frac{1}{r}=\frac{1}{r_1}+\frac{1}{r_2},$~then
\begin{align}\label{fpr}
  \||\nabla|^sG(u)\|_{L^r_x}\lesssim\|G'(u)\|_{L_x^{r_1}}\||\nabla|^su\|_{L_x^{r_2}}.
\end{align}
\end{lemma}

For the purpose of completeness, we recall some preliminaries on the Cauchy problem \eqref{maineq} as follows. (See \cite{CW} and section 2 of \cite{DM} for more details.)

\begin{lemma}\label{cwresult}

(a)[Uniqueness] Let $u$ and $\tilde{u}$ be two solutions of \eqref{maineq} on an interval $I$ containing $0$ with the same initial data $u_0$, then $u=\tilde{u}$.\vspace{2mm}
 
(b)[Existence] For initial data $u_0\in \dot{H}^1$, there exists a unique solution $u$ of \eqref{maineq} on a maximal interval $(-T_{-}(u_0),T_{+}(u_0))$. \vspace{2mm}

(c)[Finite blow-up criterion] Assume $T_{+}(u_0)<\infty$, then $||u||_{S(0,T_{+})}=+\infty$. Similar statement holds for $T_{-}(u_0)$. \vspace{2mm}

(d)[Scattering] If $T_{+}(u_0)=\infty$ and $||u||_{S(0,T_{+})}<\infty$, then $u(t)$ scatters forward in the sense that there exists $u_{+}\in \dot{H}^1$ such that
\[\lim_{t\rightarrow +\infty}||u(t)-e^{it\Delta}u_{+}||_{\dot{H}^1}=0.
\]
Similar statement holds for $T_{-}(u_0)$.\vspace{2mm}

(e)[Continuity] Let $\tilde{u}$ be a solution of \eqref{maineq} on $I$ containing $0$. Assume that for some constant $A>0$,
\[ \sup\limits_{t\in I} ||\tilde{u}||_{\dot{H}^1}+||\title{u}||_{S(I)} \leq A.
\]
Then there exists $\epsilon_0=\epsilon_0(A)>0$ and $C_0=C_0(A)$ such that for any $u_0\in \dot{H}^1$ with $||u_0-\tilde{u}_0||_{\dot{H}^1}=\epsilon<\epsilon_0$, the solution $u$ of \eqref{maineq} with initial data $u_0$ is defined on $I$ and satisfies $||u||_{S(I)} \leq C_0$ and $\sup\limits_{t\in I} ||u(t)-\tilde{u}(t)||_{\dot{H}^1} \lesssim_{C_0} \epsilon$.
\end{lemma}

\begin{lemma}[Uniform boundedness of $\dot{H}^1$-norm]\label{ubdd}
If $u$ is a subcritical threshold solution to \eqref{maineq} in the sense of \eqref{subcr}, then there exists $C>0$ such that,
\begin{equation}
 C^{-1}||u||^2_{\dot{H}^1}  \leq E(u(t))\leq C||u||^2_{\dot{H}^1}.
\end{equation}
According to the conservation law, 
\begin{equation}
||u||_{L_t^{\infty}\dot{H}_x^1(I\times \mathbb{R}^d)}<\infty,
\end{equation}
where $I$ is the maximal interval of definition. 
\end{lemma}

The proof of Lemma \ref{ubdd} is based on the following lemma (see Lemma 3.4 of \cite{KM1} for the proof of Lemma \ref{lm3.4}):
\begin{lemma}\label{lm3.4}
Let $f\in \dot{H}^1$ and $||f||_{\dot{H}^1}\leq ||W||_{\dot{H}^1}$. Then
\begin{equation}
\frac{||f||^2_{\dot{H}^1}}{||W||^2_{\dot{H}^1}} \leq \frac{E(f)}{E(W)}.
\end{equation}
In particular, $E(f)$ is positive.
\end{lemma}
\emph{Proof of Lemma \ref{ubdd}:} It is as same as Remark 2.7 of \cite{DM} so we omit it. \vspace{5mm}

As for transformations, if $v$ is a function defined on $\mathbb{R}^d$, as a convention, we write
\begin{equation}\label{symmetry}
v_{[\lambda_0,x_0]}(x)=\frac{1}{\lambda_0^{\frac{d-2}{2}}}v(\frac{x-x_0}{\lambda_0}),\textmd{ and } v_{[\theta_0,\lambda_0,x_0]}(x)=e^{i\theta_0} \frac{1}{\lambda_0^{\frac{d-2}{2}}}v(\frac{x-x_0}{\lambda_0}).
\end{equation}

Now we are ready to state the important compactness result as follows, which can be approached by a useful lemma (Lemma \ref{sequence}) which is based on profile decomposition. The difference of these results from the radial case is not big and we refer to section 2 of \cite{DM} for the radial case. 
\begin{theorem}[Compactness and global existence]\label{compactness}
Let $u$ be a subcritical solution to (1.1) in the sense of \eqref{subcr} with initial data $u(0)=u_0$ and $I$ its maximal interval of existence. Then $u$ is global, i.e. $I=\mathbb{R}$. If $S_{\mathbb{R}}(u)=\infty$, then there exists $\lambda(t)\in (0,\infty)$, $x(t) \in \mathbb{R}^d$ such that the set
\[  K:=\{\lambda(t)^{-\frac{d-2}{2}}u(t,\frac{x-x(t)}{\lambda(t)}):t\in I\}
\]
is precompact in $\dot{H}^1$. 
\end{theorem}
\emph{Remark.} Solutions that satisfy the compactness condition in Theorem \ref{compactness} are known as `almost periodic solutions'.\vspace{2mm} 

Theorem \ref{compactness} is the nonradial analogue of Theorem 2.1 of \cite{DM} and the proof consists of three steps. First, we study the properties of sequences of subcritical threshold solutions (Lemma \ref{sequence}). Second, we use the result obtained in the first step to show compactness property. At last, we use the mass concentration phenomenon to prove the global existence. Compared to the radial case, once the first step is established, then the last two steps are almost same.
\begin{lemma}\label{sequence}
Let $\{ u^0_n \}_{n\in \mathbb{N}}$ be a sequence of functions in $\dot{H}^1$ satisfying
\begin{equation}
E(u^0_n)\leq E(W),\quad ||u^0_n||_{\dot{H}^1} \leq ||W||_{\dot{H}^1}.
\end{equation}
Let $u_n$ be the solution to \eqref{maineq} with initial condition $u^0_n$. Then up to the extraction of a subsequence of $\{ u^0_n \}_{n\in \mathbb{N}}$, at least one of the following statements holds:\vspace{2mm}

(a)[Compactness] There exists sequences $x_n$ and $\lambda_n$ such that $\{(u^0_n)_{[\lambda_n,x_n]}\}_n$ converges in $\dot{H}^1$. \vspace{2mm}

(b)[Vanishing for positive direction] For every $n$, $u_n$ is defined on $[0,\infty)$ and $ \lim\limits_{n\rightarrow \infty} ||u_n||_{S([0,\infty))}=0$.\vspace{2mm}

(c)[Vanishing for negative direction] For every $n$, $u_n$ is defined on $(-\infty,0]$ and $ \lim\limits_{n\rightarrow \infty} ||u_n||_{S((-\infty,0])}=0$. \vspace{2mm}

(d)[Uniform scattering] For every $n$, $u_n$ is defined on all $\mathbb{R}$ and there is a constant $C$ independent of $n$ such that
\[ ||u_n||_{S(\mathbb{R})}  \leq C.
\]
\end{lemma}
\emph{Remark.} To overcome the gap from radial case to nonradial case, we can use a profile decomposition for nonradial case. Then the proof of Lemma \ref{sequence} follows as in Lemma 2.5 of \cite{DM} and we omit it.

\section{No high-to-low cascade and properties of the almost periodic solution}
In this section, we prove some important properties of the almost periodic solution in Theorem \ref{compactness}. First, we show that there is no high-to-low frequency cascade scenario in the sense of Theorem \ref{theorem}. Theorem \ref{theorem} is essential for us to obtain the the negative regularity of the almost periodic solution. Then we apply the negative regularity to obtain some other properties, including the $L^2$-finiteness property following the arguments in \cite{KV ecg}. 
\subsection{No high-to-low frequency cascade}
\begin{theorem}\label{theorem}
Let $u$ be as in Theorem \ref{compactness}. Thus $u$ is global. Then there exists $\lambda_0>0$ such that
\[\inf_{t\in \mathbb{R}}\lambda(t)\geq \lambda_0.\]
\end{theorem}
Theorem 3.1 is important for us to study the almost periodic solutions in Theorem \ref{compactness}. Moreover, Theorem \ref{theorem} can be obtained by using the following lemma.
\begin{lemma}\label{lemma3.2}
In Theorem \ref{theorem}, if $u\in L_t^{\infty}\dot{H}^s(\mathbb{R} \times \mathbb{R}^d)$ for some $s<1$, then the conclusion of Theorem \ref{theorem} holds.
\end{lemma}
\emph{Proof of Lemma \ref{lemma3.2}:} If not, without loss of generality, we assume that there exists $t_n\rightarrow\infty$ such that $\lambda(t_n)\rightarrow 0$. By compactness of $K$, given $\eta>0$, there exists $R=R(\eta)>0$ such that 
\begin{equation}
\int_{|\xi|\geq R\lambda(t)}|\xi|^2|\hat{u}(t,\xi)|^2d\xi<\eta.
\end{equation}
Then we have 
\begin{align*}
||u(t_n,\cdot)||_{\dot{H}^1}^2=&\int_{|\xi|\leq R\lambda(t_n)}+\int_{|\xi|\geq R\lambda(t_n)}|\xi|^2|\hat{u}(t_n, \xi)|^2d\xi\\
\leq & (R\lambda(t_n))^{2-2s}||u||_{\dot{H}^s}^2+\eta.
\end{align*}
Since $\eta>0$ is arbitrary, we obtain
$$\lim_{n\rightarrow \infty}||u(t_n,\cdot)||_{\dot{H}^1}^2=0,$$
which a contradiction. The proof of Lemma \ref{lemma3.2} is complete.\vspace{5mm}

So now it suffices to show that $u\in L_t^{\infty}\dot{H}_x^s$ for some $s<1$. This can be verified if we can show that there exists $2<p<\frac{2d}{d-2}$ such that 
\begin{equation}\label{gain}
||u(t, \cdot)||_{L^p}\lesssim 1+\lambda(t)^{-c},\quad \text{ for some sufficiently small } c>0.
\end{equation}
Let $N$ be a dyadic number. Let $a>0$ be sufficiently small which will be determined later. Let $q>\frac{2d}{d-2}$. First, we let $a$ satisfies $0<a<d(1-\frac{2}{q})-2$. Then we define
\begin{equation}
A(N):=N^{-d(\frac{1}{q}-\frac{d-2}{2d})}\sup_{t\in \mathbb{R}}\alpha(t)^{a} ||u_N(t)||_{L_x^q},
\end{equation}
where $ \alpha(t):=\min\{1, \lambda(t)\}$. We choose $q=5$ for the case $d=5$ and $q=\frac{2(d-2)}{d-4}$ for the case $d\geq 6$. Clearly, by Bernstein's inequality, we have  $A(N)\lesssim 1$.

\begin{lemma}
Let $u$ be as in Theorem \ref{compactness}. Let $\eta>0$ be small, $a>0$, then there exists $N_0(\eta,a)>0$ such that 
\begin{equation}\label{small}
\sup_{t\in \mathbb{R}}\alpha(t)^a\int_{|\xi|\leq N_0}|\xi|^2|\hat{u}(t, \xi)|^2d\xi<\eta.
\end{equation}
\end{lemma}
\emph{Proof: }By the compactness property of $K$, there exists $c_1>0$ such that 
\begin{equation}\label{low}
\sup_{t\in \mathbb{R}} \int_{|\xi|\leq c_1\lambda(t)}|\xi|^2 |\hat{u}(t,\xi)|^2d\xi<\eta.
\end{equation}
Let $\epsilon$ be a small constant to be decided.  Let $N_0$ be the largest dyadic number that is no larger than $c_1\epsilon$. Then for $\lambda(t)\geq \epsilon$, we have 
\begin{equation}
\sup_{t\in \mathbb{R}, \lambda(t)\geq \epsilon}\alpha(t)^a\int_{|\xi|\leq N_0}|\xi|^2|\hat{u}(t, \xi)|^2d\xi<\eta.
\end{equation}
For $\lambda(t)\leq \epsilon$, we have
\begin{align*}
\sup_{t\in \mathbb{R}, \lambda(t)\leq \epsilon}\alpha(t)^a\int_{|\xi|\leq N_0}|\xi|^2|\hat{u}(t, \xi)|^2d\xi\leq \epsilon^a ||W||_{\dot{H}^1}^2<\eta,
\end{align*}
 provided that we choose $\epsilon$ small enough such that $\epsilon^a ||W||_{\dot{H}^1}^2\leq\eta$.

\begin{lemma}\label{lowlow}
Let $u$ be as in Theorem \ref{compactness}, then there exists a constant $C$ such that 
\begin{equation}\label{3.7}
\lambda(t)\geq Ct^{-1/2}.
\end{equation}

\end{lemma}
 Lemma \ref{lowlow} follows as in \cite{KM1} and we omit the proof. Using the similar proof as Lemma 6.2 of \cite{KV ecg}, we have the following recurrence lemma.
\begin{lemma}
Let $\eta$ and $N_0$ be as in Lemma 3.3. For all $N\leq 10N_0$, we have 
\begin{equation}
A(N)\lesssim \Big(\frac{N}{N_0}\Big)^{\alpha}+\eta^{\frac{4}{d-2}}\sum_{\frac{N}{10}\leq N_1\leq N_0}\Big(\frac{N}{N_1}\Big)^{\alpha}A(N_1)+\eta^{\frac{4}{d-2}} \sum_{N_1<\frac{N}{10}}\Big(\frac{N_1}{N}\Big)^{\alpha}A(N_1),
\end{equation}
where $\alpha=\textmd{min }\{\frac{3d}{2}-3-a-\frac{3d}{q},\frac{3}{2}-a\}$.
\end{lemma}
\emph{Proof of Lemma 3.5:} According to time-translation symmetry, it suffices to prove 
\begin{equation}
N^{-d(\frac{1}{q}-\frac{d-2}{2d})} ||u_N(0)||_{L_x^q}\lesssim \Big(\frac{N}{N_0}\Big)^{\alpha}+\eta^{\frac{4}{d-2}}\sum_{\frac{N}{10}\leq N_1\leq N_0}\Big(\frac{N}{N_1}\Big)^{\alpha}A(N_1)+\eta^{\frac{4}{d-2}} \sum_{N_1<\frac{N}{10}}\Big(\frac{N_1}{N}\Big)^{\alpha}A(N_1).
\end{equation}
We consider the case when $d\geq 6$, the case when $d=5$ is similar. Using the following no waste Duhamel formula (see \cite{TVZ2}),
$$u(t)\xrightarrow{weakly} i\int_t^{\infty}e^{i(t-\tau)\Delta}F(u(\tau))d\tau,$$
 we have 
\begin{align*}
N^{-d(\frac{1}{q}-\frac{d-2}{2d})}||u_N(0)||_{L^q}\lesssim &N^{-d(\frac{1}{q}-\frac{d-2}{2d})} \int_0^{\infty}||e^{-it\Delta}P_N F(u(\tau))||_{L^q}d\tau\\
=&N^{-d(\frac{1}{q}-\frac{d-2}{2d})}\Big\{ \int_0^{N^{-2}}+\int_{N^{-2}}^{\infty}||e^{-it\Delta}P_N F(u(\tau))||_{L^q}d\tau\Big\}.
\end{align*}
We estimate the above two terms respectively. On one hand, by dispersive estimate,
\begin{align*}
&N^{-d(\frac{1}{q}-\frac{d-2}{2d})}\int_{N^{-2}}^{\infty}||e^{-it\Delta}P_N F(u(\tau))||_{L^q}d\tau\\
\lesssim & N^{-d(\frac{1}{q}-\frac{d-2}{2d})}\int_{N^{-2}}^{\infty}\tau^{-\frac{d}{2}(1-\frac{2}{q})}||P_N F(u(\tau))||_{L^{\frac{q}{q-1}}}d\tau\\
\lesssim &   N^{-d(\frac{1}{q}-\frac{d-2}{2d})}\int_{N^{-2}}^{\infty}\alpha(\tau)^{-a}\tau^{-\frac{d}{2}(1-\frac{2}{q})}\alpha(\tau)^{a}||P_N F(u(\tau))||_{L^{\frac{q}{q-1}}}d\tau\\
\lesssim &  N^{-d(\frac{1}{q}-\frac{d-2}{2d})}\sup_t (\alpha(t)^a ||P_NF(u(t))||_{L^{\frac{q}{q-1}}})\int_{N^{-2}}^{\infty}\tau^{\frac{a}{2}-\frac{d}{2}(1-\frac{2}{q})} d\tau\\
\lesssim & N^{-d(\frac{1}{q}-\frac{d-2}{2d})}N^{d(1-\frac{2}{q})-2-a}\sup_t (\alpha(t)^a ||P_NF(u(t))||_{L^{\frac{q}{q-1}}})\\
=& N^{\frac{3d}{2}-3-a-\frac{3d}{q}}\sup_t (\alpha(t)^a ||P_NF(u(t))||_{L^{\frac{q}{q-1}}}).
\end{align*}
On the other hand, by Bernstein's inequality,
\begin{align*}
&N^{-d(\frac{1}{q}-\frac{d-2}{2d})}\int_0^{N^{-2}} ||e^{-it\Delta}P_N F(u(\tau))||_{L^q}d\tau\\
\lesssim & N^{-d(\frac{1}{q}-\frac{d-2}{2d})+d(\frac{q-1}{q}-\frac{1}{q})}\int_0^{N^{-2}}\alpha(t)^{-a}\alpha(t)^a ||P_N F(u(\tau))||_{L^{\frac{q}{q-1}}}d\tau\\
\lesssim & N^{-d(\frac{1}{q}-\frac{d-2}{2d})+d(\frac{q-1}{q}-\frac{1}{q})}\int_0^{N^{-2}} \tau^{a/2}d\tau \cdot \sup_t (\alpha(t)^a ||P_NF(u(t))||_{L^{\frac{q}{q-1}}})\\
\lesssim &  N^{\frac{3d}{2}-3-a-\frac{3d}{q}}\sup_t (\alpha(t)^a ||P_NF(u(t))||_{L^{\frac{q}{q-1}}}).
\end{align*}
 The rest of the proof follows as in Lemma 6.2 of \cite{KV ecg} and we omit it.\vspace{5mm}

 With the above lemma, following Proposition 6.3 of \cite{KV ecg}, using Lemma 2.14 of \cite{KV ecg} and Lemma 3.4, we obtain:

\begin{theorem}\label{3.5}
 Let $u$ be as in Theorem \ref{compactness}. Then 
\begin{equation}
||u(t)||_{L_x^p}\lesssim 1+\lambda(t)^{-c}, 
\end{equation}
 for some constant $c(d,p)=a(d-2)(\frac{1}{2}-\frac{1}{p})>0$ and for $p<\frac{2d}{d-2}$ sufficiently close to $\frac{2d}{d-2}$.
\end{theorem}
At last, using Theorem \ref{3.5}, we can have
 
 \begin{theorem}\label{thm3.6}
 Let $u$ be as in Theorem \ref{compactness}. Let $d\geq 5$, then there exists $s<1$ such that $u\in L_t^{\infty}\dot{H}^s(\mathbb{R}^d)$.
\end{theorem}
\emph{Proof of Theorem \ref{thm3.6}}: Using the same argument as in \cite{KV ecg} (double Duhamel formula), we obtain
\begin{align*}
||\nabla u_N(0)||_{L_x^2}^2\lesssim & \int_0^{\infty}\int_{-\infty}^0 \min\{|t-\tau|^{-1}, N^2\}^{\frac{d}{r}-\frac{d}{2}} ||\nabla F(u(t))||_{L_x^r}||\nabla F(u(\tau))||_{L_x^r}dtd\tau,
\end{align*}
where $r=\frac{2p(d-2)}{p(d-2)+8}$. Using Theorem \ref{3.5}, Lemma \ref{chain}, and \eqref{3.7}, we obtain
\begin{equation}\label{3.11}
||\nabla F(u(t))||_{L_x^r}\lesssim 1+t^{\epsilon}, 
\end{equation}
where $\epsilon(d,p)=\frac{1}{2}a(d-2)(\frac{1}{2}-\frac{1}{p})>0$. We note that $\epsilon$ can be arbitrarily small if we choose $a$ small enough. Moreover, using \eqref{3.11}, we have
\begin{align*}
||\nabla u_N(0)||_{L_x^2}\lesssim &\int_0^{\infty}\int_{-\infty}^0 \min\{|t-\tau|^{-1}, N^2\}^{\frac{d}{r}-\frac{d}{2}} t^{\epsilon} \tau^{\epsilon} dtd\tau\\
\lesssim & \int_0^{N^{-2}}\int_{-N^{-2}}^0 N^{2(\frac{d}{r}-\frac{d}{2})} t^{\epsilon}\tau^{\epsilon}dtd\tau\\
&+ \int_0^{N^{-2}}\int_{-\infty}^{-N^{-2}}  |t-\tau|^{-(\frac{d}{r}-\frac{d}{2})} t^{\epsilon}\tau^{\epsilon}dtd\tau\\
&+ \int_{N^{-2}}^{\infty}\int_{-N^{-2}}^0 |t-\tau|^{-(\frac{d}{r}-\frac{d}{2})} t^{\epsilon}\tau^{\epsilon}dtd\tau\\
&+ \int_{N^{-2}}^{\infty}\int_{-\infty}^{-N^{-2}} |t-\tau|^{-(\frac{d}{r}-\frac{d}{2})} t^{\epsilon}\tau^{\epsilon}dtd\tau.
\end{align*}
 Note that for $t>0$, $\tau<0$, we have $|t-\tau|>t$ and $|t-\tau|>|\tau|$. Choosing $p$ sufficiently close to $\frac{2d}{d-2}$ and $a$ small enough, we can obtain
$$||\nabla u_N(0)||_{L_x^2}\lesssim N^{2s_0},$$
 for some $s_0(d,p)>0$. This completes the proof of Theorem \ref{thm3.6}. Thus the proof of Theorem \ref{theorem} is now complete as well, noticing Lemma \ref{lemma3.2}.\vspace{2mm}
 
\subsection{Properties of almost periodic solutions }We now investigate the properties of the almost periodic solution in the sense of Theorem \ref{compactness} based on Theorem \ref{theorem}. We need the mass finiteness theorem ($u \in L^{\infty}_t L^2_x$) and the control for the translation $x(t)=o(t)$ ($o(t)/t \rightarrow 0$ as $t \rightarrow \infty$) for $d \geq 5$. These results are tightly dependent on the result $\lambda(t)\geq \lambda_0$ and will be used in next section and their proofs are similar to the results in \cite{KV ecg}. We will discuss them below. The most crucial step is to obtain the negative regularity for the almost periodic solutions.

\begin{lemma}[Negative Regularity]\label{negative}
Let $u$ be a subcritical threshold solution to \eqref{maineq} and satisfying the compactness property in Theorem \ref{compactness}. Then $u\in L^{\infty}_t H^{-\epsilon}_x(\mathbb{R}\times \mathbb{R}^d)$ for some $\epsilon=\epsilon(d)>0$.
\end{lemma}

\emph{Remark.} In particular, negative regularity implies finiteness of mass according to interpolation with $\dot{H}^1$-norm (see Lemma \ref{ubdd}). The proof of Lemma \ref{negative} is tightly dependent on Theorem \ref{theorem} (no high-to-low frequency cascade) and the rest of proof follows as in Theorem 6.1 in \cite{KV ecg}, so we omit it.

\begin{lemma}[Compactness in $L^2$]\label{cmptinL2}
Let $u$ be a subcritical threshold solution to \eqref{maineq} and satisfying the compactness property in Theorem \ref{compactness}. Then for any $\eta>0$, there exists $C(\eta)>0$ such that
\begin{equation}
\sup\limits_{t\in \mathbb{R}} \int_{|x-x(t)|\geq C(\eta)} |u(t,x)|^2 dx \lesssim \eta.
\end{equation}
\end{lemma}
The proof is based on Lemma \ref{negative}, which is same as Lemma 8.3 in \cite{KV ecg}. So we omit it.

\begin{lemma}[The control of translation $x(t)$]\label{control}
Let $u$ be a subcritical threshold solution to \eqref{maineq} and satisfying the compactness property in Theorem \ref{compactness}. Then we have the following control for scaling function $x(t)$:
\begin{equation}
x(t)=o(t) \textmd{  as  } t \rightarrow \infty.
\end{equation}
\end{lemma}
\emph{Proof.} The proof is based on Lemma \ref{negative} and Lemma \ref{cmptinL2}, which is as same as Lemma 8.3 in \cite{KV ecg}. So we omit it.

\section{Exponential convergence to $W$}
 In this section, we consider a subcritical threshold solution $u$ in Theorem \ref{main}, satisfying 
\begin{equation}
E(u_0)=E(W), \quad ||u_0||_{\dot{H}^1}<||W||_{\dot{H}^1}, 
\end{equation}
and
\begin{equation}
||u||_{S(0,\infty)}=+\infty. 
\end{equation}
 The next exponential convergence theorem is very crucial for proving the main theorem  (Theorem \ref{main}). We will use it to prove the main theorem in Section 5.
\begin{theorem}\label{4main} 
 Let $u$ be a solution of (1.1) satisfying (4.1) and (4.2). Then there exist $\theta_0 \in \mathbb{R}$, $\mu_0>0$, $x_0 \in \mathbb{R}^d$ and $c,C >0$ such that 
\[ \forall t \geq 0,\quad ||u(t)-W_{[\theta_0,\mu_0,x_0]}||_{\dot{H}^1} \leq Ce^{-ct}.
\]
\end{theorem}
\begin{corollary}\label{cor4.2} There is no solution $u$ of (1.1) satisfying (4.1) and 
\begin{equation}
||u||_{S(-\infty,0)}=||u||_{S(0,\infty)}=+\infty.
\end{equation}
\end{corollary}
\emph{Remark. }Corollary \ref{cor4.2} shows that subcritical threshold solutions can not blow up in two directions.\vspace{2mm}

We define
\begin{equation}
d(f):=\left| ||f||^2_{\dot{H}^1}-||W||^2_{\dot{H}^1} \right|,
\end{equation} 
 which measures the `distance' of $f$ from $W$. The key to proving Theorem \ref{4main} is to show that
\begin{equation}\label{conv}
\lim_{t\rightarrow +\infty}d(u(t))=0. 
\end{equation}
 However, it is not easy to prove \eqref{conv} directly so we consider a weaker statement (Lemma \ref{lm4.3}) by showing the `average' of distance $d(f)$ converge to $0$ in the sense of \eqref{4.6}.
\begin{lemma}\label{lm4.3} Let $u$ be a solution of (1.1) satisfying (4.1), (4.2). Thus, $u$ is defined on $\mathbb{R}$ according to Theorem \ref{compactness}. Then 
\begin{equation}\label{4.6}
\lim_{T \rightarrow +\infty} \frac{1}{T}\int_0^T d(u(t)) dt=0.
\end{equation}
\end{lemma}
Moreover, it is obvious that Lemma \ref{lm4.3} implies:
\begin{corollary}\label{cor4.4} Under the assumptions of Lemma 4.3, there exists a sequence $t_n \rightarrow +\infty$ such that $d(u(t_n))$ converges to 0.\end{corollary}

\emph{Proof of Lemma 4.3: }Let $u$ be such a solution. By Theorem \ref{compactness}, there exists functions $\lambda(t)$ and $x(t)$ such that $K_{+}:=\{u_{[\lambda(t), x(t)]}(t),t \geq 0\}$ is relatively compact in $\dot{H}^1$.\vspace{2mm}

Let $\phi$ be a smooth, radial function such that 
\begin{equation}
\phi(r)=
\begin{cases}
r &\mbox{r $\leq$ 1} \\
0 &\mbox{r $\geq$ 2}  
\end{cases}
\end{equation}\noindent and we define $V_R(t)=\int_{\mathbb{R}^d}\psi(x)|u(x,t)|^2dx$, where $\psi(x)=R^2\phi(\frac{|x|^2}{R^2})$ for some $R >0$. Then, we have
\[\partial_t V_R(t)=4 Im \int_{\mathbb{R}^d} \phi^{'}(\frac{|x|^2}{R^2})\bar{u}x \cdot \nabla u dx.
\]
 Since $u(t,x) \in L_t^{\infty} L_x^{2} $ (according to Lemma \ref{negative}), we have 
\begin{equation}\label{4.8}
\quad |\partial_t V_R(t)| \lesssim R||\nabla u(t)||_{L^2} ||u(t)||_{L^2} \lesssim R.
\end{equation}
\begin{align*}
 \partial_{tt} V_R(t)&=4Re\int_{\mathbb{R}^d}\psi_{ij} u_i \bar{u}_jdx-\frac{4}{d}\int_{\mathbb{R}^d}(\Delta \psi)(x)|u|^{\frac{2d}{d-2}}dx-\int_{\mathbb{R}^d}(\Delta \Delta \psi)(x)|u|^2dx \\
&=\frac{16}{d-2}d(u(t))+O(\int_{|x|\geq R}(|\nabla u|^2+|u|^{\frac{2d}{d-2}})dx)+O(\int_{R\leq |x| \leq 2R}|u|^{\frac{2d}{d-2}})^{\frac{d-2}{2}}.
\end{align*}
 Here we denote $A_R(u(t))=O(\int_{|x|\geq R}(|\nabla u|^2+|u|^{\frac{2d}{d-2}})dx)+O(\int_{R\leq |x| \leq 2R}|u|^{\frac{2d}{d-2}})^{\frac{d-2}{2}}$.\vspace{2mm}

Now we want to show the following statement:\vspace{2mm}

 For $\forall \epsilon >0$, $\exists \rho_{\epsilon}>0$, such that for $\forall R >0$, $\forall t \geq 0$, $R\lambda(t) \geq \rho_{\epsilon}+x(t)$, then
\begin{equation}\label{4.9}
\partial_{tt} V_R(t) \geq \frac{16}{d-2}d(u(t))-\epsilon. 
\end{equation}
 Statement \eqref{4.9} can be justified by using the fact that $K_{+}$ is precompact in $\dot{H}^1$ and noticing that: 
\[\int_{|x|\geq R} (|\nabla u|^{2}+|u|^{\frac{2d}{d-2}})dx = \int_{|y-x(t)|\geq R\lambda(t)}(|\nabla u_{[\lambda(t),x(t)]}|^2+|u_{[\lambda(t),x(t)]}|^{\frac{2d}{d-2}})dy. 
\]
The next claim is 
\begin{equation}\label{4.10}
\lim_{t \rightarrow +\infty} t\lambda(t)-x(t) =+\infty, \end{equation}
which can be verified by invoking Theorem \ref{theorem} and Lemma \ref{control} ($x(t)=o(t)$).\vspace{2mm}

We fix $\epsilon>0$ and we use the estimates \eqref{4.8} and \eqref{4.9} with an appropriate choice of $R$. Consider the positive number $\rho_{\epsilon}$ given by \eqref{4.9}. Take $ \epsilon_{0}$ and $M_0$ such that 
\[C\epsilon_0=\epsilon,\quad M_0 \epsilon_0 =\rho_{\epsilon}
\]
 where $C$ is the corresponding constant of inequality (4.8).\vspace{2mm}

 According to \eqref{4.10}, we know that there exists $t_0$ such that for $t \geq t_0 $,
\[t\lambda(t)-x(t) \geq \rho_{\epsilon}.
\]
We consider, for $T \geq t_0$
\[ R:=\epsilon_0 T.
\]
 If $t \in [t_0,T]$, then the definitions of $R$, $M_0$ and $t_0$ imply $R\lambda(t)-x(t) \geq \epsilon_0 T\lambda(t)-x(t) \geq \rho_{\epsilon}$. Integrating \eqref{4.9} between $t_0$ and $T$  and using estimate \eqref{4.10}, we get, by the choice of $\epsilon_0$ and $R$ 
\[ \frac{16}{d-2} \int_{t_0}^{T}|d(u(t))|dt \leq CR+\epsilon(T-t_0)\leq CR+\epsilon T\leq \frac{\epsilon}{\epsilon_0} \epsilon_0 T+\epsilon T \leq 2\epsilon T.
\]
 Letting $T$ tends to $+\infty $, we obtain
\[\limsup\limits_{T \rightarrow +\infty} \frac{1}{T} \int_0^T |d(u(t))| dt \leq \frac{d-2}{8} \epsilon,
\]
 which concludes the proof of Lemma 4.3. \vspace{5mm}
 
We now work on the gap from \eqref{4.6} to \eqref{conv}. First, we introduce the orthogonal decomposition near the ground state $W(x)$. This technique was used to treat the radial case in \cite{DM}. For the nonradial case, in addition, we need to consider the partial derivative of $W(x)$, i.e. $W_j$, ($j=1,2...d$), which will appear in the orthogonal set. 
\begin{lemma} There exists $\delta_0>0$ such that for all $f$ in $\dot{H}^1$ with $E(f)=E(W)$, $d(f)<\delta_0$, there uniquely exists parameters $(\theta,\mu,x)$ in $\mathbb{R}/2\pi \mathbb{Z} \times (0,+\infty) \times \mathbb{R}^d$ with 
\[ f_{[\theta,\mu,x]}  \perp iW,W_1,\partial_j W, \quad j=1,...,d,
\]
where $W_1:=\frac{d-2}{2}W+x\cdot \nabla W$. The mapping $f \rightarrow (\theta,\mu,x)$ is $C^1$.\end{lemma}

 Let $u$ be a solution of (1.1) on an interval $I$ such that $E(u_0)=E(W)$, and on $I$, $d(u(t))< \delta_0$. According to Lemma 4.5, there exists parameter functions $\theta(t),\mu(t),x(t) $ such that 
\begin{equation}\label{ortho}
u_{[\theta(t),\mu(t),x(t)]}(t)=(1+\alpha(t))W+\tilde{u}(t)
\end{equation}
 where
\[1+\alpha(t)=\frac{1}{||W||_{\dot{H_1}}^2}(u_{[\theta(t),\mu(t),x(t)]},W)_{\dot{H}^1} , \quad \tilde{u}(t) \in \mathcal{A}^{\perp}
\]
and $\mathcal{A}:=\Big\{W, iW, W_1, \partial_j W, j=1,..., d\Big\}$. Moreover, we define $v(t)$ by 
\[ v(t):=\alpha(t)W+\tilde{u}(t)=u_{[\theta(t),\mu(t),x(t)]}(t)-W.
\]
Furthermore, we can obtain the estimates regarding the parameter functions as follows:
\begin{lemma} We consider a subcritical threshold solution $u$ defined on $I$ satisfying $d(u(t))<\delta_0$ on $I$. Taking a smaller $\delta_0$ in Theorem 4.5 if necessary, we have estimates on $I$ as follows:
\begin{equation}\label{4.6.1}
|\alpha(t)|\approx ||v(t)||_{\dot{H}^1}\approx ||\tilde{u}(t)||_{\dot{H}^1} \approx d(u(t)),
\end{equation}
\begin{equation}\label{4.6.2}
\left|\frac{x^{'}(t)}{\mu(t)}\right|+|\alpha^{'}(t)|+|\theta^{'}(t)|+\left|\frac{\mu^{'}(t)}{\mu(t)}\right|\leq C\mu(t)^2 d(u(t)).    
\end{equation}
 Also, $\alpha(t)$ and $||u(t)||_{\dot{H}^1}^2-||W||_{\dot{H}^1}^2$ have the same sign.
\end{lemma}
 The proofs of Lemma 4.5 and Lemma 4.6 will be discussed explicitly in the Appendix (Section 6). Next, we study the non-oscillatory behavior near the ground state $W(x)$, which is significant for us to obtain \eqref{conv}.
\begin{lemma}\label{timeinterval}
Let $\{t_{0n}\}_n$ and $\{t_{1n}\}_n$ be two real sequences, $\{u_{n}\}_n$ a sequence of solutions of (1.1) on $[t_{0n},t_{1n}]$ such that $u_n(t_{0n})$ satisfies assumptions (4.1) and (4.2), $\{x_n\}_n$ a sequence of functions, and $\{\lambda_n\}_n$ a sequence of positive functions such that the set: 
\[\tilde{K}=\{(u_n(t))_{[\lambda_{n}(t),x_n(t)]},n\in \mathbb{N},t\in (t_{0n},t_{1n})\}
\]
 is relatively compact in $\dot{H}^1$. Assuming 
\begin{equation}\label{4.12}
\lim\limits_{n \rightarrow +\infty} d(u_n(t_{0n}))+d(u_n(t_{1n}))=0,
\end{equation}
 then,
\begin{equation}\label{4.13}
\lim\limits_{n\rightarrow +\infty} \{\sup\limits_{t\in (t_{0n},t_{1n})} d(u(t_n))\}=0.
\end{equation}
\end{lemma}
\emph{Remark.} As for the application of Lemma \ref{timeinterval}, we often consider the following setting. Let $u$ be a solution to \eqref{maineq} satisfying (4.1) and (4.2) and parameter functions $x(t)$ and $\lambda(t)$ given by Theorem \ref{compactness}. Moreover, let $t_n$ be a sequence given by Corollary \ref{cor4.4}. Then obviously the assumptions of Lemma \ref{timeinterval} are well satisfied. Under this assumptions, if $n$ is large enough so that $d(u_n(t)) <\delta_0$ on the interval $(t_{0n},t_{1n})$, according to Lemma 4.5, we can write
\[(u_n(t))_{[\theta_n(t),\mu_n(t),x_n(t)]}(t)=(1+\alpha_n(t))W+\tilde{u}_n(t).
\]
\begin{lemma} Under the assumptions of Lemma \ref{timeinterval}, we have
\begin{equation} 
\lim\limits_{n \rightarrow +\infty} \frac{\sup\limits_{t \in (t_{0n} ,t_{1n})}\mu_n(t)}{\inf\limits_{t \in (t_{0n} ,t_{1n})}\mu_n(t)}=1.
\end{equation}
\end{lemma}
\emph{Remark.} According to the scaling invariance, it is sufficient to prove the preceding lemmas assuming 
\begin{equation}\label{assump}
\forall n, \quad \inf\limits_{t\in [t_{0n},t_{1n}]}\lambda_n(t)=1. 
\end{equation} 

 \emph{Remark.} Under the assumptions of Lemma \ref{timeinterval}, the translation functions $x_n(t)$ in the compactness argument (Theorem \ref{compactness}) and the orthogonal decomposition (Lemma 4.5) are `comparable' in the sense that the difference between them is uniformly bounded. So we will not distinguish them. Also, we will show the scalings ($\mu(t)$ and $\lambda(t)$) are also `comparable'. They will be explained in the proof of Lemma \ref{timeinterval}.\vspace{2mm}
 
 Next, we have,

\begin{lemma}\label{4.9}
If $t_n \in (t_{0n},t_{1n})$ and the sequence $\lambda_n(t_n)$ is bounded, then
\begin{equation}\label{4convto0}
 \lim\limits_{n \rightarrow +\infty} d(u_n(t_n))=0.
\end{equation}
\end{lemma} 

\begin{lemma}\label{4.10}
Let $\{u_n\}_n$ be a sequence satisfying the assumptions of Lemma 4.7 and
\[\forall n, \quad \inf\limits_{t\in [t_{0n},t_{1n}]}\lambda_n(t)=1, 
\]
 then,
 \begin{equation}\label{main410}
     \forall n,\quad \int_{t_{0n}}^{t_{1n}} d(u_n(t))dt \leq C[d(u_n(t_{0n}))+d(u_n(t_{1n}))]. 
 \end{equation}

\end{lemma}

Lemma \ref{4.10} is a key step for us to build up the exponential convergence result. Before proving Lemma \ref{4.10}, we first show that it implies the above lemmas. A brief road map for the rest of this section is as follows. First, assuming that Lemma \ref{4.10} holds, we prove Lemma \ref{4.9} and then use Lemma \ref{4.9} to prove Lemma \ref{timeinterval}. Moreover, we show Lemma 4.8. Furthermore, we give the proof of Lemma \ref{4.10}. At last, we prove Theorem \ref{4main} and Corollary \ref{cor4.2}. \vspace{2mm}

\emph{Proof of Lemma \ref{4.9}:} We may assume $1\leq \lambda_n(t_n) \leq C$, for some $C>1$. We consider $v_n(t)=u_n(t,x-x_n(t))$. So the sequence $v_n(t_n)$ is relatively compact in $\dot{H}^1$. Assuming that \eqref{4convto0} does not hold, then up to a subsequence, noticing that the distance $d(u)$ is spatial translation invariant, we have
\begin{equation}\label{eq4.1}
\lim\limits_{n \rightarrow +\infty}v_n(t_n)=v^0 \textmd{ in } \dot{H}^1({\mathbb{R}}^d),\quad d(v^0)>0,E(v^0)=E(W) \textmd{ and }||v^0||_{\dot{H}^1}<||W||_{\dot{H}^1}.
\end{equation}
Let $v$ be the solution of \eqref{maineq} with initial condition $v^0$ at time $t=0$, which is defined for $t\geq 0$. We claim that for large enough $n$, $1+t_n \leq t_{1n}$. If not, $t_{1n}\in(t_n,t_n+1)$ for an infinite number of $n$, so that extracting a subsequence, $t_{1n}-t_n$ has a limit $\tau \in [0,1]$. By the continuity of the flow of (1.1) in $\dot{H}^1$, $v_n(t_{1n})$ converges to $v(\tau)$ with $E(v(\tau))=E(W)$ and, by \eqref{4.12}, $d(v(\tau))=0$. According to Theorem \ref{sobolev}, this implies $v=W_{[\theta_0,\lambda_0,x_0]}$ for some $\theta_0$, $\lambda_0$, $x_0$, which contradicts \eqref{eq4.1}. Thus, $(t_n,t_n+1)\subset (t_{0n},t_{1n})$ holds. By \eqref{eq4.1} and the continuity of the flow of (1.1), 
\begin{equation}\label{eq4.2}
\lim\limits_{n\rightarrow +\infty} \int_{t_n}^{1+t_n} d(v_n(t)) dt=\int_0^1 d(v(t)) dt>0.
\end{equation}
However, by Lemma \ref{4.10}, $\lim\limits_{n\rightarrow +\infty} \int_{t_{0n}}^{t_{1n}} d(u_n(t)) dt=0$, which contradicts \eqref{eq4.2}. The proof of Lemma \ref{4.9} is complete.\vspace{5mm}

 \emph{Proof of Lemma 4.7:} One may assume, for every $n$, $b_n\in (t_{0n},t_{1n})$ such that
\begin{equation}\label{eq4.3}
\lim\limits_{n\rightarrow +\infty} \lambda_n(b_n)=1.
\end{equation}
According to Lemma \ref{4.9}, 
\begin{equation}
\lim\limits_{n\rightarrow +\infty} d(u_n(b_n))=0.
\end{equation}
We will use contradiction argument to show \eqref{4.13}. Without loss of generality, we assume that for some $\delta_1>0$, 
\[ \forall n, \sup\limits_{t\in (t_{0n},b_n)} d(u_n(t))\geq \delta_1>0.
\]
(The case for interval ($b_n,t_{1n}$) is similar, so we omit it). Fix $\delta_2>0$ smaller than $\delta_1$ and the constant $\delta_0$ given by Lemma 4.5. We see that there exists $a_n\in (t_{0n},b_n)$ such that 
\begin{equation}\label{eq4.23}
d(u_n(a_n))=\delta_2 \textmd{ and } \forall t\in (a_n,b_n), d(u_n(t))<\delta_2.
\end{equation}
On $(a_n,b_n)$, the modulation parameter $\mu_n$ is well defined. Moreover, noticing the relatively compactness of $\tilde{K}$ and orthogonal decomposition (to distinguish the two translation parameters, we use $x^{'}_n(t)$ for the translation parameter in the compactness argument and $x_n(t)$ for the translation parameter in the orthogonal decomposition.), the set $\cup_n \{W_{[\frac{\lambda_n(t_n)}{\mu_n(t_n)},x_n(t_n)-x^{'}_n(t_n)]}(t),t\in [a_n,b_n]\}$ must be relatively compact, which implies
\begin{equation}\label{eq4.4}
\exists C>0,\forall t\in (a_n,b_n),\quad C^{-1}\lambda_n(t) \leq \mu_n(t) \leq C\lambda_n(t) \quad \textmd{and} \quad |x_n(t)-x^{'}_n(t)|\leq C.
\end{equation}
 Using \eqref{eq4.3}, up to a subsequence, we assume that
\[ \mu_n(b_n)\rightarrow \mu_{\infty} \in (0,\infty),\textmd{ as } n\rightarrow +\infty.
\]
Now we can show by contradiction that
\begin{equation}\label{eq4.5}
\sup\limits_{n,t\in(a_n,b_n)} \mu_n(t)<\infty.
\end{equation}
If \eqref{eq4.5} does not hold, for large enough $n$, there exists $c_n\in(a_n,b_n)$ such that
\begin{equation}
\mu_n(c_n)=2\mu_{\infty},\quad \mu_n(t)<2\mu_{\infty},t\in (c_n,b_n).
\end{equation}
By Lemma \ref{4.9}, $\lim\limits_n d(u_n(c_n))=0$. Then by Lemma 4.6, we get $\left|\frac{\mu^{'}_n(t)}{\mu_n^3(t)}\right| \leq Cd(u_n(t))$. Integrating between $c_n$ and $b_n$, we get, by Lemma \ref{4.10},
\begin{equation}\label{eq4.6}
\left|\frac{1}{\mu^2_n(c_n)}-\frac{1}{\mu^2_n(b_n)}\right|\leq C\int_{c_n}^{b_{n}} d(u_n(s)) ds \rightarrow 0,\quad \textmd{as } n\rightarrow +\infty.
\end{equation}
It is a contradiction. Thus \eqref{eq4.5} holds.\vspace{2mm}

 By \eqref{eq4.5}, $\mu_n(a_n)$ is bounded. Lemma \ref{4.9} shows that $d(u_n(a_n))$ converges to $0$, contradicting \eqref{eq4.23}. The proof of Lemma 4.7 is now complete. \vspace{5mm}

\emph{Proof of Lemma 4.8:} In view of \eqref{assump} and \eqref{eq4.4}, we may assume that  
\[ \exists C>0,\forall n,\quad C^{-1}<\inf\limits_{t\in[t_{0n},t_{1n}]} \mu_n(t)<C.
\]
Since $\mu_n$ are continuous, there exists $a_n,b_n\in [t_{0n},t_{1n}]$ such that
\[\mu_{n}(a_n)=\inf\limits_{t\in[t_{0n},t_{1n}]} \mu_n(t), \quad \mu_{n}(a_n)=\sup\limits_{t\in[t_{0n},t_{1n}]} \mu_n(t).
\]
Using the bound $|\frac{\mu^{'}_n(t)}{\mu^{3}_n(t)}| \leq Cd(u_n(t))$, Lemma \ref{4.10} and Lemma \ref{timeinterval}, we obtain
\[\lim\limits_{n\rightarrow +\infty}\left|\frac{1}{\mu^2_n(a_n)}-\frac{1}{\mu^2_n(b_n)}\right|=0.
\]
Multiplying the preceding limit by $\mu^2_n(b_n)$ yields the conclusion of Lemma 4.8, noticing that $\mu_n(b_n)$ is bounded. \vspace{5mm}

Now we are ready to prove Lemma \ref{4.10} which is a key step in this section.\vspace{2mm}

\emph{Proof of Lemma \ref{4.10}:} The key elements of the proof are orthogonal decomposition and interaction Morawetz estimate. Orthogonal decomposition is very useful when we analyze functions close to the ground state $W(x)$. Interaction Morawetz estimate was first used in \cite{Iteam1} by J. Colliander, M. Keel, G. Staffilani, H. Takaoka and T. Tao. Later, there are some modified versions of interaction Morawetz estimate used in many papers (see \cite{Iteam2,BD3,BD ecg} as examples). In particular, we apply the version of the interaction Morawetz estimate used in \cite{BD ecg} by B. Dodson. \vspace{2mm} 

We define a function $\psi \in C_0^{\infty}(\mathbb{R})$, $\psi$ even, $\psi(x)=1$ for $|x| \leq 1$ and $\psi(x)=0$ for $|x| > 2$. We let 
\[ \phi(x-y)=\int \psi^2(|x-s|) \psi^2(|y-s|) ds.
\]
We note that $\phi$ is supported on $|x|\leq 4$ and we define the interaction Morawetz action as follows
\[ M_{R,n}(t)=\int |u_n(t,y)|^2 \phi(\frac{x-y}{R})(x-y)_j \cdot  Im[\bar{u}_n\partial_j u_n](t,x) dxdy.\]

The main idea of proving Lemma \ref{4.10} is to prove an upper bound for $M_{R,n}(t)$ and a lower bound for $\partial_t M_{R,n}(t)$. Then we can integrate $\partial_t M_{R,n}(t)$ to obtain the conclusion \eqref{main410}.\vspace{2mm}

\emph{Step 1:} a bound from above for $M_{R,n}$. In this step, we show that there exists some constant $C>0$ such that
\begin{equation}\label{star1}
\forall R >0 ,\forall n, \forall t \in (t_{0n},t_{1n}),\quad |M_{R,n}(t)|\leq CR^2 d(u_n(t)).
\end{equation} 
When $d(u_n(t))$ is big, \eqref{star1} can be verified by showing $|M_{R,n}(t)|\leq CR^2 ||u_n(t)||_{\dot{H}^1}$ which can be proved by using H$\ddot{o}$lder's inequality. When $d(u_n(t))$ is small, we apply the orthogonal decomposition to write $(u_n(t))_{[\theta_n(t),\mu_n(t),x_n(t)]}=W+v_n(t)$, with $||v_n(t)||_{\dot{H}^1} \leq Cd(u_n(t))$. Using the change of variable $x=\frac{x^{'}-x_n(t)}{\mu_n(t)}$, $y=\frac{y^{'}-x_n(t)}{\mu_n(t)}$,
\begin{align*}
M_{R,n}(t)&=\mu_n(t)^{-2d} \int |u_n(t,\frac{y^{'}-x_n(t)}{\mu_n(t)})|^2 \phi(\frac{x^{'}-y^{'}}{R\mu_n(t)}) (\frac{x^{'}-y^{'}}{\mu_n(t)})\cdot Im[\bar{u}_{n} \partial_j u_n](t,\frac{x^{'}-x_n(t)}{\mu_n(t)})d x^{'}d y^{'}\\
&=R^2 \mu_n(t) \int \mu_n(t)^{-d} |u_n(t,\frac{y^{'}-x_n(t)}{\mu_n(t)})|^2 \frac{1}{(R\mu_n(t))^2} \phi(\frac{x^{'}-y^{'}}{R\mu_n(t)}) (\frac{x^{'}-y^{'}}{\mu_n(t)})\\ 
&\cdot Im[\frac{1}{\mu_n(t)^{\frac{d-2}{2}}} \bar{u}_n \frac{1}{\mu_n(t)^{\frac{d}{2}}} \partial_j u_n](t,\frac{x^{'}-x_n(t)}{\mu_n(t)})d x^{'}d y^{'}.
\end{align*}
For the quantity above, we write:
\[ Im[(W+\bar{v}_n)\nabla (W+v_n)]=Im(W\nabla v_n+\bar{v}_n \nabla W+\bar{v}_n \nabla v_n).\]
 And by using the boundedness of $\mu_n(t)$, the mass finiteness of $u_n(t)$ and Cauchy-Schwarz inequality, we get the bound $|M_{R,n}(t)|\leq CR^2(|v_n(t)|_{\dot{H}^1}+|v_n(t)|_{\dot{H}^1}^2)$ which yields \eqref{star1}, for $d(u_n(t))\leq \delta_1$, $\delta_1$ small. Then the proof of \eqref{star1} is complete.\vspace{2mm}

 Step 2: a bound from below for $\partial_t M_{R,n}(t)$. We want to show that there exists some constant $C^{'}$,
\begin{equation}\label{star2}
\exists R_0, \forall R \geq R_0,\forall n,\forall t \in (t_{0n},t_{1n}),\quad  \partial_t M_{R,n}(t) \geq C^{'}d(u_n(t)).
\end{equation}
We can use \eqref{star1} and \eqref{star2} to prove Lemma \ref{4.10}. Indeed, integrating \eqref{star2} between $t_{0n}$ and $t_{1n}$ we get
\[ C^{'}\int_{t_{0n}}^{t_{1n}} d(u_n(t))dt\leq M_{R_0,n}(t_{0n})+M_{R_0,n}(t_{1n})
\]
which implies Lemma \ref{4.10} in view of \eqref{star1}. Thus, it suffices to prove \eqref{star2}. First, we calculate $\partial_t M_{R,n}(t)$. For convenience, we use $M_{R}(t)$ instead of $M_{R,n}(t)$ by considering $u(t)$ instead of $u_n(t$). The estimates will work for all $u_n$ with same constants since they are solutions to the same initial value problem \eqref{maineq}. A direct calculation shows that
\begin{equation}\label{Mora}
\aligned
\partial_t M_{R}(t)&=(d-2)\int \phi(\frac{x-y}{R})  |u(t,y)|^2[|\nabla u(t,x)|^{2}-|u(t,x)|^{\frac{2d}{d-2}}]dxdy\\
&- (d-2)\int Im[\bar{u}\partial_j u](t,y) \phi(\frac{x-y}{R}) Im[\bar{u}\partial_j u](t,x) dxdy \\
&+2\int (\phi^{'}(\frac{x-y}{R})\frac{(x-y)_k (x-y)_j}{R|x-y|})[Re(\partial_j \bar{u}\partial_k u)(t,x)-\frac{1}{d}\delta_{jk}|u(t,x)|^{\frac{2d}{d-2}}]|u(t,y)|^2 dxdy\\
&-2\int (\phi^{'}(\frac{x-y}{R})\frac{(x-y)_k (x-y)_j}{R|x-y|})Im(\bar{u}\partial_j u)(t,x)Im(\bar{u}\partial_k u)(t,y)dxdy\\
&-\frac{1}{2}\int (\Delta[d\phi(\frac{x-y}{R})+\phi^{'}(\frac{x-y}{R})\frac{|x-y|}{R}]) |u(t,x)|^2 |u(t,y)|^2dxdy \\
&=A+B+C+D+E.
\endaligned
\end{equation}
For the first term we can take $d(u_n(t))$ out by noticing 
\[\int |\nabla u|^2-|u|^{2^*} =\frac{2}{d-2}d(u(t)).
\]
Furthermore, we can use the $L^2$-finiteness of $u_n$ to write $\partial_t M_{R,n}(t)$ to be the sum of a main term and a remainder term as follows,
\begin{equation}\label{4.28}
\partial_t M_{R,n}(t)=C^{''}d(u_n(t))+A_R(u_n(t)),
\end{equation}
where $C^{''}$ depends on the dimension $d$ and the initial data $u_0$ in \eqref{maineq}.\vspace{2mm} 

The next step is to control the remainder term $A_R(u_n(t))$ depending on the distance $d(u_n(t))$. We claim the following two estimates for the remainder term $A_R(u_n(t))$:
\begin{equation}\label{est1}
\forall \epsilon>0,\quad \exists \rho_{\epsilon}>0 \textmd{ such that } \forall n, \forall t \in (t_{0n},t_{1n}), \forall R\geq \frac{\rho_{\epsilon}}{\lambda_n(t)},\quad |A_R(u_n(t))| \leq \epsilon.    
\end{equation}
\begin{equation}\label{est2}
\aligned
&\exists \delta_2,\quad \forall n, \forall t \in (t_{0n},t_{1n}),\forall R\geq \frac{1}{\mu_n(t)}, \\ 
&d(u_n(t))\leq \delta_2 \textmd{  implies  }|A_R(u_n(t))| \lesssim  d(u_n(t))+d(u_n(t))^2,
\endaligned
\end{equation}
where $c$ is a positive constant.\vspace{2mm} 

We prove \eqref{est1} first. For any $\epsilon>0$, we can choose $R$ big enough to ensure that the terms $C,D,E$ in \eqref{Mora} are arbitrarily small. For term $B$, we can write it as 
\begin{align*}
B=&- (d-2)\int Im[\bar{u}\partial_j u](t,y) Im[\bar{u}\partial_j u](t,x) dxdy\\
&- (d-2)\int Im[\bar{u}\partial_j u](t,y) (\phi(\frac{x-y}{R})-1) Im[\bar{u}\partial_j u](t,x) dxdy\\
:=&B_1+B_2.
\end{align*}
$B_2$ is obvious small when $R$ is big enough. Note that 
$$B_1=-(d-2)P(u(t))^2=-(d-2)(\sum_{j=1}^d P_j(u(0)))^2,$$
where $P_j(u):=\int Im[\bar{u}\partial_j u]$ is the $j$-th component of the momentum. Because $u\in L^2(\mathbb{R}^d)$, by conservation of momentum, we have $P_j(u(t))=P_j(u(0))$.  As a consequence, unless $P_j(u)=0$, we cannot expect $\tilde{B}_1$ to have any decay in time. So we need to use $A$ to beat $B$. We write $A$ as 
\begin{align*}
A=& (d-2)\int  |u(t,y)|^2[|\nabla u(t,x)|^{2}-|u(t,x)|^{\frac{2d}{d-2}}]dxdy\\
&+(d-2)\int (\phi(\frac{x-y}{R})-1)  |u(t,y)|^2[|\nabla u(t,x)|^{2}-|u(t,x)|^{\frac{2d}{d-2}}]dxdy\\
:=& A_1+A_2.
\end{align*}
It's obvious that $A_2$ is small when $R$ is big enough. We consider the quantity 
$$F(u(t))=\int |\nabla u(t,x)|^2 |u(t,y)|^2-Im[\bar{u}\partial_j u](t,x)Im[\bar{u}\partial_j u](t,y)dxdy.$$
It's straightforward to see that for any given $\xi=(\xi_1, ...,\xi_d)\in \mathbb{R}^d$, we have  $F(e^{ix\cdot \xi}u(t))=F(u(t))$.  Note that 
$$\int Im[e^{-ix\cdot\xi}\bar{u}\partial_j (e^{ix\cdot \xi}u)]=\xi_j ||u(t)||_{L^2}^2+\int Im[\bar{u}\partial_j u].$$
Since $||u_0||_{L^2}^2\neq 0$, we choose 
$$\xi_j:=-\frac{\int Im[\bar{u}\partial_j u]}{||u(t)||_{L^2}^2}.$$
\noindent So we have 
$$\int Im[e^{-ix\cdot\xi}\bar{u}\partial_j (e^{ix\cdot \xi}u)]=0.$$
 By conservation of momentum and mass, $\xi_j$ is independent of $t$. So we have 
\begin{align*}
F(u(t))=\int |\nabla (e^{ix\cdot\xi}u(t,x))|^2 |u(t,y)|^2dxdy.
\end{align*}
And
\begin{align*}
 A_1+B_1=&(d-2)F-(d-2)\int |u(t,y)|^2 |u(t,x)|^{\frac{2d}{d-2}}dxdy\\
=& (d-2)\int [ |\nabla (e^{ix\cdot\xi}u(t,x))|^2-|u(t,x)|^{\frac{2d}{d-2}}] |u(t,y)|^2 dxdy. 
\end{align*}
 Recall the following lemma (which is a consequence of Lemma 3.4 of \cite{KM1}),
\begin{lemma}\label{trap}
\noindent Assume that 
$$||\nabla u||_{L^2}^2\leq (1-\delta)||\nabla W||_{L^2}^2,$$
\noindent where $\delta>0$. Then there exists $\bar{\delta}=\bar{\delta}(\delta, d)$ such that 
$$\int |\nabla u|^2-|u|^{\frac{2d}{d-2}}\geq \bar{\delta} \int |\nabla u|^2.$$
\end{lemma}
Then we have,
\begin{corollary}\label{cor4.12}
 Assume 
\begin{equation}\label{ass1}
\int |\nabla u|^2dx\leq (1-\delta)||\nabla W||_{L^2}^2,
\end{equation} 
then
$$\int |\nabla (e^{ix\cdot\xi}u(t,x))|^2-|u(t,x)|^{\frac{2d}{d-2}}dx\geq \bar{\delta} \int|\nabla u|^2.$$
\end{corollary}
\emph{Proof of Corollary \ref{cor4.12}:} If $\int |\nabla (e^{ix\cdot \xi}u(t,x))|^2 dx<\int |\nabla u|^2dx$, then 
$$\int |\nabla (e^{ix\cdot \xi}u(t,x))|^2<(1-\delta)||\nabla W||_{L^2}^2.$$
 By Lemma \ref{trap}, we have 
\begin{equation}
\int |\nabla (e^{ix\cdot\xi}u(t,x))|^2-|u(t,x)|^{\frac{2d}{d-2}}dx\geq \bar{\delta} \int |\nabla( e^{ix\cdot \xi}u(t,x))|^2 dx.
\end{equation}
If $\int |\nabla (e^{ix\cdot \xi}u(t,x))|^2 dx\geq \int |\nabla u|^2dx$, then 
\begin{align*}
&\int |\nabla (e^{ix\cdot\xi}u(t,x))|^2-|u(t,x)|^{\frac{2d}{d-2}}dx\geq \int |\nabla u|^2-|u|^{\frac{2d}{d-2}}dx\\
=&\frac{2}{d-2}d(u(t))=\frac{2}{d-2}\Big| ||\nabla u||_{L^2}^2-||\nabla W||_{L^2}^2\Big|\geq \delta ||\nabla W||_{L^2}^2\geq \bar{\delta}||\nabla u||_{L^2}^2,
\end{align*}
 where the last inequality is by (\ref{ass1}). The proof of Corollary \ref{cor4.12} is complete.\vspace{5mm}

Using Corollary \ref{cor4.12}, we have 
\begin{equation}
A_1+B_1\geq \bar{\delta}||u_0||_{L^2}^2 ||\nabla u||_{L^2}^2\geq C^{''} d(u(t)).
\end{equation}
Thus the proof of \eqref{est1} is complete.\vspace{5mm}

Now we turn to \eqref{est2}. According to Lemma 4.5, when the distance $d(u(t))$ is small enough, we decompose $u$ as $$u(t)_{[\theta(t),\mu(t),x(t)]}=W+v(t),$$ \noindent with $||v(t)||_{\dot{H}^1} \leq Cd(u(t))$. We use the change of variables $x=\frac{x^{'}-x(t)}{\mu(t)}$, $y=\frac{y^{'}-x(t)}{\mu(t)}$. In new variables, we write $A$ as
\begin{align*}
A=&(d-2)\mu(t)^{d-2}\int\phi(\frac{x'-y'}{R\mu(t)})    |\frac{1}{\mu(t)^{(d-2)/2}}u(t,\frac{y'-x(t)}{\mu(t)})|^2\\
&\times  [|\nabla u_{[\theta,\mu,x]}(t,x')|^2-|u_{[\theta,\mu,x]}(t,x')|^{\frac{2d}{d-2}}]dx'd \frac{y'}{\mu}\\
=&\frac{d-2}{\mu(t)^2}\int \phi(\frac{x'-y'}{R\mu(t)})|u_{[\theta,\mu,x]}(t,y')|^2 [|\nabla u_{[\theta,\mu,x]}(t,x')|^2-|u_{[\theta,\mu,x]}(t,x')|^{\frac{2d}{d-2}}]dx'dy'\\
=& \frac{d-2}{\mu(t)^2}\int \phi(\frac{x'-y'}{R\mu(t)})|(W+v(t,y')|^2 [|\nabla (W+v)(t,x')|^2-|(W+v)(t,x')|^{\frac{2d}{d-2}}]dx'dy'\\
:=& A(W+v).
\end{align*}
We denote $A(W)$ by replacing $W+v$ by $W$.  We do the same for the terms $B, C, D, E$. We have 
\begin{align*}
B=&-\frac{d-2}{\mu(t)^2}\int Im[\overline{(W+v)}\partial_j (W+v)](t,y') \phi(\frac{x'-y'}{R\mu(t)}Im[\overline{W+v}\partial_j(W+v)](t,x')dx'dy'\\
:=&B(W+v).
\end{align*}

\begin{align*}
C=&\frac{2}{\mu(t)^2}\int (\phi'(\frac{x'-y'}{R\mu(t)})\frac{(x'-y')_k (x'-y')_j}{R|x'-y'|\mu(t)})[Re(\partial_j \overline{W+v}\partial_k (W+v))(t,x')\\
& -\frac{1}{d}\delta_{jk}|(W+v)(t,x)|^{\frac{2d}{d-2}}]|(W+v)(t,y)|^2 dx'dy'\\
:=& C(W+v).
\end{align*}

\begin{align*}
D=&-\frac{2}{\mu(t)^2}\int (\phi'(\frac{x'-y'}{R\mu(t)})\frac{(x'-y')_k(x'-y')_j}{R|x'-y'|\mu(t)})\\
&\times   Im[\overline{(W+v)}\partial_j (W+v)](t,y)Im[\overline{W+v}\partial_j (W+v)](t,x')dx'dy'\\
:=& D(W+v).
\end{align*}

\begin{align*}
E=&-\frac{1}{2\mu(t)^2}\int (\Delta [d\phi(\frac{x'-y'}{R\mu(t)})+\phi'(\frac{x'-y'}{R\mu(t)})\frac{|x'-y'|}{R\mu(t)}])|(W+v)(t,y)|^2 |(W+v)(t,x)|^2 dx'dy'\\
:=&E(W+v).
\end{align*}
We also write $M_R(t)$ as $M_R(W+v)$, and we denote $M_R(W)$ by replacing $W+v$ by $W$. Note that 
$$A(W)+B(W)+C(W)+D(W)+E(W)=0.$$
So we have 
\begin{align*}
\frac{d}{dt}M_R(t)=&(A(W+v)-A(W))+(B(W+v)-B(W))+(C(W+v)-C(W))\\
&+(D(W+v)-D(W))+(E(W+v)-E(W)).
\end{align*}
We will estimate the above terms respectively.

\subsection{Estimate $A(W+v)-A(W)$.} Note that 
\begin{align*}
A(W+v)=& \frac{d-2}{\mu(t)^2}\int|(W+v)(t,y')|^2 [|\nabla (W+v)(t,x')|^2-|(W+v)(t,x')|^{\frac{2d}{d-2}}]dx'dy'\\
&+ \frac{d-2}{\mu(t)^2}\int (\phi(\frac{x'-y'}{R\mu(t)})-1)|(W+v)(t,y')|^2 [|\nabla (W+v)(t,x')|^2-|(W+v)(t,x')|^{\frac{2d}{d-2}}]dx'dy'\\
:=& A_1(W+v)+A_2(W+v).
\end{align*}
Using the fact that 
$$\int |\nabla W|^2-|W|^{\frac{2d}{d-2}}=d(W)=0,$$
we have
\begin{align*}
A(W)=& \frac{d-2}{\mu(t)^2}\int |W(t,y')|^2 [|\nabla W(x')|^2-|W(t,x')|^{\frac{2d}{d-2}}dx'dy'\\
&+\frac{d-2}{\mu(t)^2} \int  (\phi(\frac{x'-y'}{R\mu(t)})-1)|W(t,y')|^2 [|\nabla W(x')|^2-|W(t,x')|^{\frac{2d}{d-2}}]dx'dy'\\
:=& A_1(W)+A_2(W).
\end{align*}
 We have $A_1(W)=0$. Making change of variables back, we have 
\begin{align*}
A_1(W+v)=& (d-2)\int |u(t,y)|^2 [|\nabla u(t,x)|^2-|u(t,x)|^{\frac{2d}{d-2}}]dx'dy'\\
=& (d-2)\int |u(t,y)|^2dy d(u(t))\\
\geq (d-2)&||u_0||_{L^2}^2 d(u(t)).
\end{align*}
 For $A_2(W+v)-A_2(W)$, exploring the cancellations, we have 
\begin{align*}
&A_2(W+v)-A_2(W)\\
=& \frac{d-2}{\mu(t)^2}\int  (\phi(\frac{x'-y'}{R\mu(t)})-1)\Big\{|(W+v)(t,y')|^2 [|\nabla (W+v)(t,x')|^2-|(W+v)(t,x')|^{\frac{2d}{d-2}}]\\
&-|W(y')|^2 [|\nabla W(x')|^2-|W(t,x')|^{\frac{2d}{d-2}}]\Big\} dx'dy'\\
=&  \frac{d-2}{\mu(t)^2}\int  (\phi(\frac{x'-y'}{R\mu(t)})-1)\Big\{(|(W+v)|^2-|W|^2) [|\nabla (W+v)(t,x')|^2-|(W+v)(t,x')|^{\frac{2d}{d-2}}] dx'dy'\\
&+ \frac{d-2}{\mu(t)^2}\int  (\phi(\frac{x'-y'}{R\mu(t)})-1)|W(y')|^2 \Big\{[|\nabla (W+v)(t,x')|^2-|(W+v)(t,x')|^{\frac{2d}{d-2}}]\\
 &- [|\nabla W(x')|^2-|W(t,x')|^{\frac{2d}{d-2}}]\Big\} dx'dy'\\
:=& \tilde{A}_1+\tilde{A}_2.
\end{align*}
 Note that 
$$|W+v|^2-|W|^2=2Re[W\bar{v}]+|v|^2.$$
 Also, if we denote 
$$\mu_0:=\inf_{t\in \mathbb{R}}\mu(t),$$
 then
$$supp(\phi(\frac{x'-y'}{R\mu(t)})-1)\subset \Big\{ (x', y')\quad | \quad |x'|\geq \frac{R\mu_0}{2}, or~|y'|\geq \frac{R\mu_0}{2}\Big\}.$$
 Note that the energy center of $W+v$ is zero. Given $\epsilon$ sufficiently small, we can choose $R$ sufficiently large such that
\begin{align*}
\int_{|x|\geq R\mu_0/2}|W(x)|^2dx\leq \epsilon^2,
\end{align*}
and
\begin{align*}
\int_{|x|\geq R\mu_0/2}|\nabla (W+v)|^2+|W+v|^{\frac{2d}{d-2}}dx\leq \epsilon^2.
\end{align*}
 So we have 
\begin{align*}
\tilde{A}_1\leq & 2\frac{d-2}{\mu_0^2}\int_{|x|\geq R\mu_0/2}[|vW(t,y')|+|v(t,y')|^2dy'\\
&\times  \int|\nabla (W+v)(t,x')|^2+|(W+v)(t,x')|^{\frac{2d}{d-2}}dx'\\
+& 2\frac{d-2}{\mu_0^2}\int [|vW(t,y')|+|v(t,y')|^2dy'\\
&\times  \int_{|x|\geq R\mu_0/2}|\nabla (W+v)(t,x')|^2+|(W+v)(t,x')|^{\frac{2d}{d-2}}dx'\\
\leq &  c||W||_{L^2(|x|\geq R\mu_0/2)}||v||_{L^2(|x|\geq R\mu_0/2)}+c||v||_{L^2}^2\\
\leq & c\epsilon ||v||_{\dot{H}^1}+c||v||_{\dot{H}^1}^2\\
\leq & c\epsilon d(u(t))+cd(u(t))^2.
\end{align*}
The estimate for $\tilde{A}_2$ is similar and we omit it. So we obtain
\begin{equation}\label{AA}
A(W+v)-A(W)\geq  \frac{||u_0||_{L^2}^2}{2}d(u(t))-c\epsilon d(u(t))-cd(u(t))^2.
\end{equation}

 \subsection{Estimate $B(W+v)-B(W)$.} Since $W$ is real, we have 
\begin{equation}
\begin{split}
& Im[\overline{(W+v)}\partial_j  (W+v)]\\
=&Im[\overline{v}\partial_j W]+Im[\overline{v}\partial_j v]+Im[\overline{W}\partial_j v].
\end{split}
\end{equation}

 Integrating by parts, if the derivative is taken on $W$, we move the derivative to take on $v$, and obtain
\begin{equation}\label{BB}
\begin{split}
|B(W+v)-B(W)|\leq & C||\nabla v||_{L^2}^2(||W||_{L^2}^2+||W||_{L^2}||v||_{L^2}+||v||_{L^2}^2)\\
\leq & Cd(u(t))^2.
\end{split}
\end{equation}

 \subsection{Estimate for $C(W+v)-C(W)$.} Similar to the case of $A_2(W+v)-A_2(W)$, we can obtain
\begin{equation}\label{DD}
|C(W+v)-C(W)\leq C\epsilon d(u(t))+Cd(u(t))^2.
\end{equation}

 \subsection{Estimate for $D(W+v)-D(W)$.} This is similar to that of $B(W+v)-B(W)$. We have 
\begin{equation}\label{EE}
|D(W+v)-D(W)|\leq Cd(u(t))^2.
\end{equation}

 \subsection{Estimate for $E(W+v)-E(W)$.} This can be estimated in a way similar to the previous cases. The worst part in this case is 
\begin{align*}
E_1:=\int (\Delta [d\phi(\frac{x'-y'}{R\mu(t)})+\phi'(\frac{x'-y'}{R\mu(t)})\frac{|x'-y'|}{R\mu(t)}])|(W+v)(t,y)|^2 |v(t,x)|^2 dx'dy',
\end{align*}
 because we don't know whether $||v||_{L^2}\lesssim d(u(t))$. This is not a problem, because we have 
$$(\Delta [d\phi(\frac{x'-y'}{R\mu(t)})+\phi'(\frac{x'-y'}{R\mu(t)})\frac{|x'-y'|}{R\mu(t)}])\lesssim \frac{1}{1+|x'-y'|^2}.$$
 By Hardy's inequality, we obtain
$$E_1\lesssim Cd(u(t))^2.$$
 A careful analysis gives 
\begin{equation}\label{FF}
|E(W+v)-E(W)|\lesssim \epsilon d(u(t))+d(u(t))^2.
\end{equation}
 To sum up, by (\ref{AA}), (\ref{BB}), (\ref{DD}), (\ref{EE}), (\ref{FF}), we obtain
\begin{equation}
\frac{d}{dt}M_R(t)\geq (d+2)||u_0||_{L^2}^2 d(u(t))-\epsilon  d(u(t))-d(u(t))^2.
\end{equation}
Then take $R$ sufficiently large so take $\epsilon$ is sufficiently small such that $\epsilon \ll \frac{1}{2}||u_0||_{L^2}^2$. We obtain
\begin{equation}
\frac{d}{dt}M_R(t)\geq \frac{1}{4}||u_0||_{L^2}^2 d(u(t))-d(u(t))^2.
\end{equation}
This completes the estimate for \eqref{star2}.\vspace{2mm}

At last, we use the estimates \eqref{est1} and \eqref{est2} to prove \eqref{star2}. According to \eqref{est2}, there exists some $\delta_3>0,R_1>0$ such that for $d(u_n(t))\leq \delta_3,R\geq R_1$,
\begin{equation}
|A_{R}(u_n(t))| \leq \frac{C^{''}}{2} d(u_n(t)).
\end{equation}
Now we use \eqref{est1} with $\epsilon=\frac{C^{''}\delta_3}{2}$, we obtain $|A_{R}(u_n(t))| \leq \frac{C^{''}}{2} d(u_n(t))$ for $d(u_n(t)) \geq \delta_3,R\geq R_2$. Estimate \eqref{star2} holds with $R_0:=\textmd{max}\{R_1,R_2\}$ and $C^{'}=\frac{C^{''}}{2}$ in view of \eqref{4.28}. The proof of Lemma \ref{4.10} is complete.\vspace{5mm}

Now we will show Theorem 4.1 and Corollary 4.2. The proof of Theorem 4.1 consists of three steps. We consider $u$ be a solution of (1.1) satisfying (4.1) and (4.2).\vspace{2mm}

 \emph{Step 1 (convergence of $d(u(t))$ to $0$): }First, we prove (4.5). By Corollary 4.4, there exists a strictly increasing sequence $\{t_n\}_{n\in \mathbb{N}}$ such that: 
\[ \lim\limits_{n\rightarrow +\infty} t_n =+\infty, \quad \lim\limits_{n\rightarrow +\infty} d(u(t_n))=0. 
\]
We let $t_{0n}=t_n$, $t_{1n}=t_{n+1}$, and $\lambda_n(t)=\lambda(t)$, where $\lambda(t)$ is given by Theorem \ref{compactness} (compactness argument). Then the assumptions of Lemma \ref{timeinterval} are satisfied by the sequences $\{u_n\}_{n\in \mathbb{N}}$, $\{t_{0n}\}_{n\in \mathbb{N}}$, $\{t_{1n}\}_{n\in \mathbb{N}}$ and $\{\lambda_n\}_{n\in \mathbb{N}}$. Thus,
\[\lim\limits_{n\rightarrow +\infty}(\sup\limits_{t\in [t_n,t_{n+1}]} d(u(t)))=0
,\]
 which implies (4.5).\vspace{2mm}

According to Lemma 4.5 (orthogonal decomposition), we decompose $u$ for large $t$ as follows,
\[u_{[\theta(t),\mu(t),x(t)]}(t)=(1+\alpha(t))W+\tilde{u}(t),\quad \tilde{u}(t)\in \mathcal{A}^{\perp}.
\]
 The conclusion of Theorem 4.1 is equivalent to the existence of $\mu_{\infty}>0$, $\theta_{\infty}\in \mathbb{R}$, $x_{\infty}\in \mathbb{R}^{d}$ and $C,c>0$ such that 
\begin{equation}\label{eq4.7}
d(u(t))+|\alpha(t)|+||\tilde{u}(t)||_{\dot{H}^1}+|x(t)-x_{\infty}|+|\theta(t)-\theta_{\infty}|+|\mu(t)-\mu_{\infty}| \leq Ce^{-ct}.
\end{equation}

\emph{Step 2 (convergence of $\mu(t)$) :} In this step, we show by contradiction that $\mu(t)$ has a limit $\mu_{\infty}\in (0,+\infty)$ as $t\rightarrow +\infty$. This step is essential for us to control other parameter functions. If not, log$(\mu(t))$ does not satisfy the Cauchy criterion as $t\rightarrow +\infty$, which implies that there exists two sequences $\{T_n\},\{T^{'}_n\} \rightarrow +\infty$ such that
\begin{equation}\label{4.30}
\lim\limits_{n\rightarrow +\infty} \frac{|\mu(T_n)|}{|\mu(T^{'}_n)|}=L\neq 1.
\end{equation}
Without loss of generality, we assume $T_n<T^{'}_{n}$. The Step 1 shows that $d(u(T_n))$ and $d(u(T^{'}_n))$ tend to $0$. Now we let $u_n=u$, $t_{0n}=T_n$, $t_{1n}=T^{'}_n$ and $\lambda_n(t)=\lambda(t)$, where $\lambda(t)$ is again given by Theorem \ref{compactness}. Then the assumptions of Lemma 4.8 are satisfied, which shows
\[ \lim\limits_{n\rightarrow +\infty} \frac{\textmd{inf}_{T_n \leq t \leq T^{'}_n} \mu(t)}{\textmd{sup}_{T_n \leq t \leq T^{'}_n} \mu(t)}=1.
\]
This contradicts \eqref{4.30}. Thus
\begin{equation}
\lim\limits_{t\rightarrow +\infty} \mu(t)=\mu_{\infty}\in (0,\infty).
\end{equation}

 \emph{Step 3 (Proof of Theorem 4.1) :} We are now ready to prove \eqref{eq4.7}. First, we show that $d(u(t))$ converges exponentially to $0$. We claim the following inequality
\begin{equation}\label{eq4.8}
\exists C>0,\forall t \geq 0 , \int_{t}^{+\infty} d(u(\tau))d\tau\leq Cd(u(t)).
\end{equation}
 If \eqref{eq4.8} does not hold, there exists a sequence $T_n \rightarrow +\infty$ such that 
\begin{equation}\label{eq4.9}
\int_{T_n}^{+\infty} d(u(\tau))d\tau\geq nd(u(T_n)).
\end{equation}
As shown in Step 2, $\mu(t)$ is bounded from below. This implies that the parameter $\lambda(t)$ of Theorem \ref{compactness} is also bounded from below. By Step 1, the assumptions of Lemma \ref{4.10} are satisfied for the sequence $\{u_k\}_{k}$, with $k=(n,n^{'})$, and $u_k=u(t)$, $\lambda_k(t)=\lambda(t)$, $t_{0k}=T_n$ and $t_{1k}=T_{n^{'}}$. Thus
\[\forall n,n^{'},n<n^{'},\quad \int_{T_{n}}^{T_{n^{'}}} d(u(t)) dt \leq C[d(u(T_n))+d(u(T_{n^{'}}))] .
\]
 We see $\int_{T_n}^{+\infty} d(u(t))dt\leq Cd(T_n)$, which contradicts \eqref{eq4.9}. So we know \eqref{eq4.8} holds.\vspace{2mm}

Now by \eqref{eq4.8} we have, for some constants $C,c>0$
\[\int_{t}^{+\infty}d(u(\tau)d\tau \leq Ce^{-ct}.
\]
Together with the estimate $|\alpha^{'}(t)| \leq Cd(u(t))$ of Lemma 4.6, we obtain
\[|\alpha(t)|=\left|\int_t^{+\infty}\alpha^{'}(\tau) d\tau \right|\leq Ce^{-ct}.
\]
By Lemma 4.6, we know $|\alpha^{'}(t)| \approx d(u(t))$ which gives us the bound on $d(u(t))$ in \eqref{eq4.7}. Moreover, using Lemma 4.6 again, we can obtain the bounds on $|\alpha^{'}(t)|$ and $||\tilde{u}(t)||_{\dot{H}^1}$ in \eqref{eq4.7}. Now it is left to show the exponential convergence of $\theta(t)$, $\mu(t)$ and $x(t)$ in \eqref{eq4.7}.\vspace{2mm}

Actually, it suffices to prove the exponential convergence for $\theta^{'}(t)$, $\mu^{'}(t)$ and $x^{'}(t)$ respectively according to fundamental theorem of Calculus and Cauchy criterion. Eventually, using the estimate \eqref{4.6.2} $\left|\frac{x^{'}(t)}{\mu(t)}\right|+|\alpha^{'}(t)|+|\theta^{'}(t)|+\left|\frac{\mu^{'}(t)}{\mu(t)}\right|\leq C\mu(t)^2 d(u(t))$ of Lemma 4.6 and the boundedness of $\mu(t)$, we can obtain \eqref{eq4.7}. Now the proof of Theorem 4.1 is complete.\vspace{5mm}

At last, we prove Corollary 4.2.\vspace{2mm}

\emph{Proof of Corollary 4.2:} It suffices to show that there is no solution $u$ of (1.1) satisfying (4.1) and (4.3). Let $u$ be such a solution. By applying Theorem 4.1 forward and backward, the set $\{u(t),t\in \mathbb{R}\}$ is relatively compact in $\dot{H}^1$. Moreover, we have
\[\lim\limits_{t\rightarrow +\infty}d(u(t))=\lim\limits_{t\rightarrow -\infty}d(u(t))=0.
\]
According to Lemma \ref{4.10} with $u_n(t)=u(t)$, $t_{0n}=-n$, $t_{1n}=n$ and $\lambda_n(t)=1$, we have $\int_{-\infty}^{+\infty}d(u(t))dt=\lim\limits_{n\rightarrow +\infty}\int_{-n}^{n}d(u(t))dt=0$. This implies $d(u_0)=0$, which clearly contradicts (4.1). The proof of Corollary 4.2 is complete.

\section{Proof of main result}
In this section, we prove the main theorem of this paper, i.e. Theorem \ref{main}. The following proposition will be applied. 

\begin{proposition}\label{5.1}
Let $C,c>0$. Assume $u$ is the solution of \eqref{maineq} satisfying $E(u)=E(W)$, $||u_0||_{\dot{H}^1}<||W||_{\dot{H}^1}$ and 
\begin{equation}
||u(t)-W||_{\dot{H}^1} \leq Ce^{-ct},\quad \forall t\geq 0.
\end{equation}
Then there exists $T\in \mathbb{R}$ such that
\begin{equation}
u(t)=W^{-}(t+T).
\end{equation}
\end{proposition}

\emph{Remark.} The proof of Proposition \ref{5.1} can be found in Lemma 6.5 (and Corollary 6.6) of \cite{DM} for $d=5$ and Theorem 4.1 (and Corollary 4.2) of \cite{LZ} for $d\geq 6$. Regarding the construction of $W^-$, we refer to section $6$ of \cite{DM} for $d=5$ and section $4$ of \cite{LZ} for $d\geq 6$. \vspace{2mm}

Now we prove the main theorem of this paper based on Proposition \ref{5.1} and Theorem \ref{4main} as follows.\vspace{2mm}

\emph{Proof of Theorem \ref{main}:} Let $u$ be a maximal-lifespan solution to \eqref{maineq} on $I$ satisfying $E(u)=E(W)$ and $||u_0||_{\dot{H}^1}<||W||_{\dot{H}^1}$. Then according to Theorem \ref{compactness}, we have $I=\mathbb{R}$. If $u$ scatters in both time directions, that is one case. If not, without loss of generality, assuming $u$ blows up forward in time, using Theorem \ref{4main}, we conclude that there exists $\theta_0,\mu_0,x_0$, $C>0,c>0$ such that  
\begin{equation}
||u(t)-W_{[\theta_0,\mu_0,x_0]}||_{\dot{H}^1} \leq Ce^{-ct},\quad \forall t\geq 0.
\end{equation}
This implies
\begin{equation}
||u_{[-\theta_0,\mu_0^{-1},-x_0]}(t)-W||_{\dot{H}^1} \leq Ce^{-c\mu_0^{-2}t}.
\end{equation}
By Proposition \ref{5.1}, we conclude that there exists $T\in \mathbb{R}$ such that $u_{[-\theta_0,\mu_0^{-1},-x_0]}=W^{-}(t+T)$.\vspace{2mm}

Thus, we get $u(t,x)=e^{i\theta_0}\mu_0^{-\frac{d-2}{2}}W^{-}(\mu_0^{-2}t+T,\mu_0^{-1}(x-x_0))$.\vspace{2mm}

This shows that $u=W^{-}$ up to symmetries. The proof of the main theorem is complete.

\section{Appendix}

In this section, we give proofs for Lemma 4.5 and Lemma 4.6 in Section 4. We refer to Lemma 3.6 and Lemma 3.7 of \cite{DM} for the analogues of these two lemmas in the radial setting. The main idea of the proofs are similar to the radial case and the difference is to deal with and control the translation parameter function $x(t)$. First, we use implicit function theorem to prove the orthogonal decomposition near ground state $W(x)$ based on the properties of radial functions. And then we prove the estimates regarding parameter functions.

\begin{lemma}\label{odd}
Let $f\in L^1(\mathbb{R}^d)$ be a radial function. Then $\int \frac{x_i}{|x|}f(x)dx=0$.
\end{lemma}
\emph{Proof:} This follows from the observation that $\frac{x_i}{|x|}f(x)$ is an odd function in $x_i$-direction. 

\begin{lemma}
 The set 
$$\mathcal{A}:=\Big\{W, iW, W_1, \partial_j W, j=1,..., d\Big\}$$
\noindent is an orthogonal set in $\dot{H}^1(\mathbb{R}^d)$.
\end{lemma}
\emph{Proof of Lemma 6.2:} It's clear that $\partial_j W$ is orthogonal to $iW$. Note that $W, W_1$ are radial. So $\Delta W, \Delta W_1$ are also radial. Then $\partial_j W \Delta W$ is odd in $x_j$ direction. So we have 
\begin{align*}
(W, \partial_j W)_{\dot{H}^1}=\int \partial_j W \Delta W =0.
\end{align*}
 Similarly, we have 
$(W_1, \partial_j W)_{\dot{H}^1}=0$ and $ (\partial_j W, \partial_k W)_{\dot{H}^1}=0, ~for~j\neq k.$\vspace{5mm}

Now we give the proof of Lemma 4.5 as follows:\vspace{2mm}

\emph{Proof of Lemma 4.5:} The proof is almost the same as the proof of Lemma 3.6 in \cite{KM1} (radial case). Define for $x=(x_1,...., x_d)$,
$$J_0:(\theta, \mu, x_1,...,x_d,f)\mapsto (f_{[\theta, \mu, x]}, iW)_{\dot{H}^1}, \quad J_1:(\theta, \mu, x_1,...,x_d,f)\mapsto (f_{[\theta, \mu, x]}, W_1)_{\dot{H}^1},$$
$$G_j: (\theta, \mu, x_1,...,x_d,f)\mapsto (f_{[\theta, \mu, x]}, \partial_j W)_{\dot{H}^1},\quad j=1,...,d.$$
 We have 
$$\frac{\partial J_0}{\partial \theta}(0,1,0,W)=\int |\nabla W|^2, \quad \frac{\partial{J_0}}{\partial\mu}(0, 1,0, W)=0,\quad \frac{\partial J_0}{\partial x_j}(0,1,0,W)=0,$$
$$\frac{\partial J_1}{\partial \theta}(0,1,0,W)=0, \quad\quad \frac{\partial{J_1}}{\partial\mu}(0, 1, 0,W)=-\int |\nabla W_1|^2, \quad \frac{\partial J_1}{\partial x_j}(0,1,0,W)=0,$$
$$\frac{\partial G_k}{\partial \theta}(0,1,0,W)=0, \quad\quad \frac{\partial{G_k}}{\partial\mu}(0, 1, 0,W)=0, \quad \frac{\partial G_k}{\partial x_j}(0,1,0,W)=-\delta_{jk}\int |\partial_j W|^2.$$
Also, we have 
\[J_0(0,1,0,W)=J_1(0,1,0,W)=G_j(0,1,0,W)=0,\quad j=1,...,d.\]

By implicit function theorem, there exists $\epsilon_0, \eta_0>0$ such that for $h\in \dot{H}^1$, if $||h-W||_{\dot{H}^1}<\epsilon_0$, then there exists unique $(\theta, \mu, x)$,
$$|\theta|+|\mu-1|+|x|\leq \eta_0,$$
 and
$$(h_{[\theta, \mu, x]}, iW)_{\dot{H}^1}=(h_{[\theta, \mu, x]}, W_1)_{\dot{H}^1}=(h_{[\theta, \mu, x]},\partial_j W)_{\dot{H}^1}=0.$$
Thus, the proof of Lemma 4.5 is complete.\vspace{5mm} 
 
 As a result, given $f\in \dot{H}^1$ such that $E(f)=E(W)$, by variational characterization of $W$, if $d(f)$ is sufficiently small, there exist parameters $(\theta, \mu, x)$ and $g$ such that 
$$f_{[\theta, \mu, x]}=W+g.$$ 
Moreover, we consider $u$ be a subcritical threshold solution to (1.1) on $I$ such that $d(u(t))<\delta_0$ where $\delta_0$ is given in Lemma 4.5. Then there exists time-dependent parameter functions $\alpha(t), \mu(t), x(t)$ and $\tilde{u}(t)\in \mathcal{A}^{\perp}$ such that
$$u_{[\theta(t), \mu(t), x(t)]}=(1+\alpha(t))W+\tilde{u}(t),$$
where we denote 
$$v(t)=\alpha(t)W+\tilde{u}(t).$$
 
 Additionally, we have estimates for the parameter functions. We give proof of Lemma 4.6 as follows.\vspace{2mm}

 \emph{Proof of Lemma 4.6:} The proof of \eqref{4.6.1} is almost same as the radial case and we refer to section 7 of \cite{DM} for the details. For simplicity, we recall and use the notation $U:=u_{[\theta(t), \mu(t), x(t)]}$ which means 
$$U(t,x)=\frac{e^{i\theta(t)}}{\mu(t)^{(d-2)/2}}u(t, \frac{x-x(t)}{\mu(t)}).$$
We apply change of variables $t=t(s)$ such that 
$$\frac{dt}{ds}=\frac{1}{\mu(t)^2}.$$
In order to prove \eqref{4.6.2}, it suffices to prove
\begin{equation}\label{goal}
\left|\frac{x_s}{\mu}(s)\right|+|\alpha_s(s)|+|\theta_s(s)|+\left|\frac{\mu_s}{\mu}(s)\right|\leq Cd(u(s)).    
\end{equation}
By straight calculations, we have 
\begin{align*}
\partial_s U(t,x)=&e^{i\theta}\Big\{-\frac{d-2}{2}\frac{\mu_s}{\mu^{d/2}}u(t,\frac{x-x(t)}{\mu})+\frac{1}{\mu^{(d-2)/2}}\frac{dt}{ds} u_t(t,\frac{x-x(t)}{\mu})\\&+\frac{1}{\mu^{(d-2)/2}}\partial_s (\frac{x-x(t)}{\mu(t)})\cdot \nabla u(t,\frac{x-x(t)}{\mu})\Big\}+i\theta_s U\\
=&-\frac{d-2}{2}\frac{\mu_s}{\mu}U(t,x)+\frac{e^{i\theta}}{\mu^{(d+2)/2}}u_t(t, \frac{x-x(t)}{\mu})-e^{i\theta}\frac{\mu_s}{\mu^{(d+2)/2}}(x-x(t))\cdot \nabla u(t,\frac{x-x(t)}{\mu})+i\theta_sU\\
&-e^{i\theta}\frac{\partial_s x(t(s))}{\mu^{d/2}}\cdot \nabla u(t, \frac{x-x(t)}{\mu})\\
=&-\frac{\mu_s}{\mu}\Big(\frac{d-2}{2}U(t,x)+(x-x(t))\cdot \nabla U(t,x)\Big)+\frac{e^{i\theta}}{\mu^{(d+2)/2}}u_t(t, \frac{x-x(t)}{\mu})-\frac{x_s}{\mu}\cdot \nabla U(t,x)+i\theta_s U.
\end{align*}
Also, we have
$$\Delta U(t,x)=\frac{e^{i\theta}}{\mu^{(d+2)/2}}\Delta u(t,\frac{x-x(t)}{\mu}).$$
So we obtain
\begin{equation}\label{tildeu}
\begin{split}
&i\partial_s U+\Delta U\\
=&-i\frac{\mu_s}{\mu}\Big(\frac{d-2}{2} U(t,x)+(x-x(t))\cdot \nabla U(t,x)\Big)+\frac{e^{i\theta}}{\mu^{(d+2)/2}}(iu_t+\Delta u)(t,\frac{x-x(t)}{\mu})-i\frac{x_s}{\mu}\cdot \nabla U(t,x)-\theta_s U\\
=&-i\frac{\mu_s}{\mu}\Big(\frac{d-2}{2}U(t,x)+(x-x(t))\cdot \nabla U(t,x)\Big)+\frac{e^{i\theta}}{\mu^{(d+2)/2}}|u|^{\frac{4}{d-2}}u(t,\frac{x-x(t)}{\mu})-i\frac{x_s}{\mu}\cdot \nabla U(t,x)-\theta_s U \\
=& -i\frac{\mu_s}{\mu}\Big(\frac{d-2}{2}U(t,x)+(x-x(t))\cdot \nabla U(t,x)\Big)+|U|^{\frac{4}{d-2}}U-i\frac{x_s}{\mu}\cdot \nabla U(t,x)-\theta_s U.
\end{split}
\end{equation}
Now we decompose $U$ as 
\begin{equation}
U=W+\alpha(t)W+\tilde{u}:=W+v,\quad\quad \tilde{u}:=g_1+ig_2\in \mathcal{A}^{\perp},
\end{equation} 
where $g_1,g_2$ are real. Similar as the radial case, we rewrite (\ref{tildeu}) as
\begin{equation}\label{newtildeu}
\begin{split}
\partial_s v+&\mathcal{L}v+R(v)-\theta_s iW-i\theta_s v-i\theta_s v+\frac{\mu_s}{\mu}(W_1-x(t)\cdot \nabla W)\\&+\frac{\mu_s}{\mu}(\frac{d-2}{2}v+(x-x(t))\cdot \nabla v)-\frac{x_s}{\mu}\cdot \nabla W-\frac{x_s}{\mu}\cdot \nabla v=0,
\end{split}
\end{equation}
 where $\mathcal{L}$ and $\mathcal{R}$ are defined in section 5 and section 7 of \cite{DM}. We denote $p_c:=\frac{d+2}{d-2}$ and \eqref{newtildeu} can be further written as
\begin{equation}\label{assa}
\begin{split}
&\partial_s g_1+i\partial_s g_2+\alpha_s W+(\Delta+W^{p_c-1})g_2-i(\Delta+p_cW^{p_c-1})g_1-i\alpha(p_c-1)W^{p_c} \\
&-\theta_s iW+\frac{\mu_s}{\mu}W_1-(\frac{\mu_s}{\mu}x(t)+\frac{x_s}{\mu})\cdot\nabla W\\
=&-R(v)+i\theta_s v-\frac{\mu_s}{\mu}\Big(\frac{d-2}{2}v+x\cdot \nabla v\Big)+(\frac{\mu_s}{\mu}x(t)+\frac{x_s}{\mu})\cdot \nabla v:=\mathcal{E}.
\end{split}
\end{equation}
 Denote
$$\epsilon(s):=|d|(|d|+|\theta_s(s)|+\left|\frac{\mu_s}{\mu}(s)\right|+\left|\frac{\mu_s}{\mu}x(t)+\frac{x_s}{\mu}\right|),\quad d:=d(u(s))$$
and
$$c=||W||_{\dot{H}^1}^2,\quad \quad c_1=||W_1||_{\dot{H}^1}^2.$$
We multiply by $\Delta W$ on both sides of \eqref{assa}, integrate and then take real parts, using the facts that 
\begin{equation}\label{6ortho}
(g_1, W)_{\dot{H}^1}=(W, W_1)_{\dot{H}^1}=(W, \partial_j W)_{\dot{H}^1}=0,    
\end{equation}
we obtain
\begin{equation}\label{alphas}
c\alpha_s=-(\Delta g_2, W)_{\dot{H}^1}-(W^{p_c-1}g_2, W)_{\dot{H}^1}+O(\epsilon(s)).
\end{equation}
Multiplying by $\Delta iW$ on both sides of \eqref{assa}, integrating and taking imaginary part, we obtain
\begin{equation}\label{thetas}
c\theta_s=-(\Delta g_1, W)_{\dot{H}^1}-p_c(W^{p_c-1}g_1, W)_{\dot{H}^1}-\alpha (p_c-1)(W^{p_c}, W)_{\dot{H}^1}+O(\epsilon(s)).
\end{equation}
Multiplying by $\Delta W_1$ on both sides of (\ref{assa}), integrating and then taking real part, we obtain
\begin{equation}\label{mus}
c_1\frac{\mu_s}{\mu}=-(\Delta g_2, W_1)_{\dot{H}^1}-(W^{p_c-1}g_2, W_1)_{\dot{H}^1}+O(\epsilon(s)).
\end{equation}
Multiplying by $\Delta\partial_j W$ on both sides of (\ref{assa}), integrating and then taking real part, denoting 
$$\lambda_j:=||\partial_j W||_{\dot{H}^1}^2, \quad \quad \beta_j:=\frac{\mu_s}{\mu}x_j(t)+\frac{\partial_s x_j(t)}{\mu},$$
and recalling \eqref{6ortho}, we obtain
\begin{equation}\label{betajs}
\lambda_j \beta_j(s)=-((\Delta+W^{p_c-1})g_2, \partial_j W)_{\dot{H}^1}+O(\epsilon(s)).
\end{equation}
Putting (\ref{alphas}), (\ref{thetas}), (\ref{mus}), (\ref{betajs}) together, we obtain
\begin{equation}\label{6.11}
|\alpha_s|+|\theta_s|+\left|\frac{\mu_s}{\mu}\right|+\sum_j|\beta_j(s)|\leq C(||g||_{\dot{H}^1}+O(\epsilon(s)))\leq Cd(u)+O(\epsilon(s)).
\end{equation}
Let $\delta_0$ be small enough, we obtain
\begin{equation}\label{6.12}
|\alpha_s|+|\theta_s|+\left|\frac{\mu_s}{\mu}\right|+\sum_j|\beta_j(s)|\leq Cd(u),
\end{equation}
which implies \eqref{goal}.\vspace{2mm} 

In particular, as shown in the proof of Theorem 4.1, we can obtain exponential estimate for $d(u)$ and the convergence of $\mu$. Based on these, according to \eqref{6.12}, we can obtain the exponential estimates for all the parameter functions as follows
\begin{equation}\label{6.13}
|\alpha_s|+|\theta_s|+|\mu_s|+d(u(t)) \leq Ce^{-\delta t}.    
\end{equation}  
Also, noticing that $x(t)\lesssim t$, in \eqref{6.11}, we can obtain
\begin{equation}\label{6.14}
\left|\frac{x_s}{\mu}\right|\lesssim d(u(t))+\left|\frac{\mu_s}{\mu}x(t)\right| \lesssim e^{-\delta t}+te^{-\delta t}\lesssim e^{-\frac{1}{2}\delta t}.  
\end{equation} 
The proof of Lemma 6.6 is complete.\vspace{5mm}

\noindent \textbf{Acknowledgments.} The authors would like to express thanks to Professor Benjamin Dodson and Chenjie Fan for useful discussions, suggestions and comments.

\hfill \linebreak
\noindent \author{Qingtang Su}

\noindent \address{Department of Mathematics, University of Michigan} \\
{530 Church Street, Ann Arbor, MI 48109-1043, U.S.}

\noindent \email{qingtang@umich.edu}\vspace{5mm}

\noindent \author{Zehua Zhao}

\noindent \address{Department of Mathematics, Johns Hopkins University}\\
{3400 N. Charles Street, Baltimore, MD 21218, U.S.}

\noindent \email{zzhao25@jhu.edu}\\

\end{document}